\documentclass{article}
\usepackage{amssymb,amsmath}

\newlength{\centrage}
\divide\centrage by 2
\oddsidemargin \centrage

\newtheorem{theorem}{Theorem}
\newtheorem{lemma}{Lemma}
\newtheorem{proposition}{Proposition}
\newtheorem{remark}{Remark}
\newtheorem{example}{Example}
\newtheorem{examples}[example]{Examples}
\newcommand{\proof}{\noindent {\sc Proof. }}
\newcommand{\qed}{\hfill\ensuremath{\Box}\vskip\baselineskip}

\newcommand{\boite}[1]{\parbox{.3\linewidth}{\begin{center}\mbox{}\\[-2.5ex]#1\end{center}}}
\newcommand{\boitelarge}[1]{\parbox{.33\linewidth}{\begin{center}\mbox{}\\[-2.5ex]#1\end{center}}}

\title{\bf A trick around Fibonacci, Lucas and Chebyshev}
\author{\small\sc Aim\'e LACHAL\footnote{Email: aime.lachal@insa-lyon.fr
\hspace{1em}URL: http://maths.insa-lyon.fr/$\mbox{}^{\sim}$lachal}\\
\small\it Universit\'e de Lyon/Institut Camille Jordan%
\footnote{Postal address: Institut National des Sciences Appliqu\'ees de Lyon,
P\^ole de Math\'ematiques, B\^atiment L\'eonard de Vinci,
20 avenue Albert Einstein, 69621 Villeurbanne Cedex, France}}
\date{}

\begin{document}
\maketitle

\begin{abstract}
In this article, we present a trick around Fibonacci numbers which can be found
in several magic books. It consists in computing quickly the sum of the successive
terms of a Fibonacci-like sequence. We give explanations and extensions of this trick to
more general sequences. This study leads us to interesting connections between
Fibonacci, Lucas sequences and Chebyshev polynomials.
\end{abstract}

\section{Introduction}

We describe a famous magic trick which can be found in many magic treatises;
see, e.g., \cite{fulves}, \cite{gardner}, \cite{larayne}, \cite{simon}.
A magician asks you to choose two integers (not too complicated
since you will have to compute several sums) in your head.
Add them and write the result for instance on a paper that the magician
does not see for the time being.
Next add the number you just wrote to the second number you have chosen and
write the result on your paper.
Next add the two numbers you wrote and write the result,
and so on. Each time, you add the two last numbers you obtained and
you write the result.
Repeat this procedure until you have written 8 numbers on your paper.
Finally add all these numbers including the two first numbers you chose:
this is the sum of 10 numbers. Now ask the magician to look at your paper.
Without knowing the two first integers you chose, he will quickly obtain the
sum by using only one term written on the paper.
Indeed, we claim that this is 11 times the fifth term you wrote
(or the fourth term starting from the last one you wrote).
You have certainly understood that you were dealing with the famous Fibonacci
numbers (with general initial conditions)...

The explanation is the following one: you choose two integers $x_0$
and $x_1$ (actually the trick holds for any complex numbers) and you
successively write
$x_2=x_0+x_1$, $x_3=x_1+x_2$, $x_4=x_2+x_3$, $x_5=x_3+x_4$,
$x_6=x_4+x_5$, $x_7=x_5+x_6$, $x_8=x_6+x_7$ and $x_9=x_7+x_8$.
All these numbers can be expressed by means of $x_0$ and $x_1$:
\begin{eqnarray*}
&x_2=x_0+x_1,\; x_3=x_0+2x_1,\; x_4=2x_0+3x_1,\; x_5=3x_0+5x_1,&
\\
&x_6=5x_0+8x_1,\; x_7=8x_0+13x_1,\; x_8=13x_0+21x_1,\; x_9=21x_0+34x_1.&
\end{eqnarray*}
Finally, you compute the sum of all the numbers from $x_0$ to $x_9$:
\begin{align*}
S=x_0+x_1+x_2+x_3+x_4+x_5+x_6+x_7+x_8+x_9.
\end{align*}
You can easily check that $S=55x_0+88x_1=11x_6$.
A simpler trick consists in writing only 4 numbers ($x_2,x_3,x_4,x_5$)
and computing the sum of 6 numbers:
\begin{align*}
S'=x_0+x_1+x_2+x_3+x_4+x_5.
\end{align*}
This time, we have $S'=8x_0+12x_1=4x_4$. We claim that the sum $S'$
is 4 times the third term you wrote.

The relationships $S=11x_6$ and $S'=4x_4$, which hold for any choice of
$x_0$ and $x_1$, are the starting point of our investigation. In this article,
we address two problems:
\begin{itemize}
\item
Is it possible to find in a Fibonacci sequence with general initial conditions
other couple of integers
$(n,m)$ such that the sum of the $n$ first terms is proportional to the
$m$th term, possibly through a simple ratio (e.g., an integer)?
\item
Is it possible to observe a similar property for more general sequences
satisfying a second order linear recurrence relation?
\end{itemize}

Fibonacci numbers consist without any doubt of one of the most well-known and important
sequence of numbers in Mathematics, Life and Nature! It enjoys so many remarkable,
fascinating, funny properties that it is the source of a lot of entertainment.
Our problem relies on certain worthwhile properties concerning the sum of Fibonacci numbers.

Although starting from a classical trick used in magic, our aim is not to
invent more spectacular magic tricks but to introduce well-known or less-known
interesting relationships between different and important mathematical objects.

Throughout the paper, $\mathbb{N}$ will denote the set of non-negative integers,
that is $\mathbb{N}=\{0,1,2,\dots\}$, $\mathbb{R}$ the set of real numbers and
$\mathbb{C}$ that of complex numbers.

The paper is organized as follows.
In Section~\ref{section-Fibo}, we recall and exhibit several features of
classical Fibonacci and Lucas numbers (satisfying the recurrence $x_{n+2}=x_{n+1}
+x_n$) and we give a positive answer to the first question we addressed above.
In Section~\ref{section-Lucas}, we consider more general Fibonacci-like
sequences, those satisfying the second order linear recurrence $x_{n+2}=
ax_{n+1}+bx_n$. For those corresponding to the value $b=-1$, we are able to
answer positively to the second question. Our analysis involves interesting
connections between such sequences and the famous Chebyshev polynomials
which are described in Section~\ref{section-cheb}.
Concerning the value $b=1$ (and $a\neq 1$), we are able to answer partially
to the second question. More precisely, we get a proportionality identity
only for certain numbers $a$.
This case is related to the family of the Fibonacci polynomials
which are also described in Section~\ref{section-cheb}.
In Section~\ref{section-conclusion}, we sum up the results obtained in this article.
Finally, we include in an appendix some additional and technical informations concerning
the case where $b=1$.

In order to facilitate the reading of the article and to make it self-contained
and broadly accessible, we report many details. Nonetheless, we assume that
the reader is familiar with linear recursive sequences.
Since there is a huge literature on Fibonacci, we have chosen to give a
concise and uniform list of references based on electronic
resources like Wikipedia and Wolfram.

\section{The classical Fibonacci and Lucas sequences}\label{section-Fibo}

\subsection{Recall of the definitions and certain properties}

In this section, we recall definitions and some well-known facts on
the classical Fibonacci and Lucas numbers; we refer the reader to
\cite{fibonb}, \cite{lucasnb}, \cite{wikifibonb} and~\cite{wikilucasnb}.
They are defined by
\begin{align*}
\mathcal{F}_0=0,\,\mathcal{F}_1=1 \;\text{ and }\;
\mathcal{F}_{n+2}=\mathcal{F}_{n+1}+\mathcal{F}_n \;\text{ for } n\in\mathbb{N},
\\
\mathcal{L}_0=2,\,\mathcal{L}_1=1 \;\text{ and }\;
\mathcal{L}_{n+2}=\mathcal{L}_{n+1}+\mathcal{L}_n \;\text{ for } n\in\mathbb{N}.
\end{align*}
We write the first terms:
\begin{eqnarray*}
&\mathcal{F}_0=0,\;\,\mathcal{F}_1=1,\; \mathcal{F}_2=1,\; \mathcal{F}_3=2,\; \mathcal{F}_4=3,
\mathcal{F}_5=5,\; \mathcal{F}_6=8,\; \mathcal{F}_7=13,\; \mathcal{F}_8=21,\; \mathcal{F}_9=34\dots&
\\
&\mathcal{L}_0=2,\;\,\mathcal{L}_1=1,\; \mathcal{L}_2=3,\; \mathcal{L}_3=4,\; \mathcal{L}_4=7,
\mathcal{L}_5=11,\; \mathcal{L}_6=18,\; \mathcal{L}_7=29,\; \mathcal{L}_8=47,\; \mathcal{L}_9=76\dots&
\end{eqnarray*}
Let us introduce the famous golden number $\varphi_1=\big(1+\sqrt5\,\big)/2$
together with
$\varphi_2=\big(1-\sqrt5\,\big)/2$. The numbers $\varphi_1$ and $\varphi_2$ are the solutions
of the quadratic equation $\varphi^2=\varphi+1$. They satisfy $\varphi_1+\varphi_2=1$
and $\varphi_1\varphi_2=-1$.
Binet's formula supplies explicit expressions for $\mathcal{F}_n$ and $\mathcal{L}_n$:
\begin{align*}
\mathcal{F}_n=\frac{\varphi_1^n-\varphi_2^n}{\varphi_1-\varphi_2}
\;\text{ and }\; \mathcal{L}_n=\varphi_1^n+\varphi_2^n.
\end{align*}

We first state a property concerning the sum of the Fibonacci numbers.
Set $\mathcal{S}_0=0$ and $\mathcal{S}_n=\sum_{k=0}^{n-1} \mathcal{F}_k$ for any
integer $n\ge 1$. The sum $\mathcal{S}_n$ can be immediately evaluated
by performing a telescopic sum: indeed, we have
$\mathcal{S}_n=\sum_{k=0}^{n-1} (\mathcal{F}_{k+2}-\mathcal{F}_{k+1})$,
which can be simplified into
\begin{align*}
\mathcal{S}_n=\mathcal{F}_{n+1}-1.
\end{align*}

In relation to the trick described in the introduction, we mention
an interesting relationship
between Fibonacci and Lucas numbers: for any $p,q\in\mathbb{N}$ such that $p\ge q$,
\begin{align}\label{rel-fibo-pq}
\mathcal{F}_{p+q}+(-1)^q\mathcal{F}_{p-q}=\mathcal{F}_p\mathcal{L}_q.
\end{align}
The proof of this identity is quite elementary:
\begin{align*}
(\varphi_1-\varphi_2)\mathcal{F}_p\mathcal{L}_q
& =(\varphi_1^p-\varphi_2^p)(\varphi_1^q+\varphi_2^q)
=(\varphi_1^{p+q}-\varphi_2^{p+q})
+(\varphi_1^p\varphi_2^q-\varphi_1^q\varphi_2^p)
\\
& =(\varphi_1^{p+q}-\varphi_2^{p+q})
+(\varphi_1\varphi_2)^q(\varphi_1^{p-q}-\varphi_2^{p-q})
=(\varphi_1^{p+q}-\varphi_2^{p+q})+(-1)^q(\varphi_1^{p-q}-\varphi_2^{p-q})
\\
& =(\varphi_1-\varphi_2)(\mathcal{F}_{p+q}+(-1)^q\mathcal{F}_{p-q}).
\end{align*}
Let $m\in\mathbb{N}$. In particular, (\ref{rel-fibo-pq}) yields for $p=q=m$ and
for $p=m+1$, $q=m$, $m$ being odd in the second case:
\begin{align}\label{rel-fibolucas}
\mathcal{F}_{2m}=\mathcal{F}_m\mathcal{L}_m\;\text{ and }\;
\mathcal{F}_{2m+1}=\mathcal{F}_{m+1}\mathcal{L}_m+1,
\end{align}
from which we derive
\begin{align*}
\mathcal{S}_{2m-1}=\mathcal{F}_m\mathcal{L}_m-1
\;\text{ and }\;\mathcal{S}_{2m}=\mathcal{F}_{m+1}\mathcal{L}_m.
\end{align*}
We then state the following important identity for our problem.
%
\begin{theorem}
If $n$ is a multiple of\/ $4$ plus $2$, we have the following equality:
\begin{align}\label{relS-fibo}
\mathcal{S}_n=\mathcal{L}_{n/2}\mathcal{F}_{n/2+1}.
\end{align}
\end{theorem}
%
We retrieve our introductory examples: $\mathcal{S}_6=\mathcal{L}_3\mathcal{F}_4
=4\mathcal{F}_4$ and $\mathcal{S}_{10}=\mathcal{L}_5\mathcal{F}_{6}=11\mathcal{F}_6$.

The identity~(\ref{relS-fibo}) can be extended to a sum starting from $\mathcal{F}_r$
for a certain $r\in\mathbb{N}$: with the same assumption on $n$, (\ref{rel-fibo-pq})
yields for $p=r+n/2+1$ and $q=n/2$:
\begin{align*}
\sum_{k=0}^{n-1}\mathcal{F}_{k+r}=\mathcal{F}_{n+r+1}-\mathcal{F}_{r+1}
=\mathcal{L}_{n/2}\mathcal{F}_{r+n/2+1}.
\end{align*}
For instance, we can find in~\cite{wikifibonb} the proportionality identities:
\begin{align*}
\mathcal{F}_{r}+\mathcal{F}_{r+1}+\dots+\mathcal{F}_{r+5} & =4\mathcal{F}_{r+4},
\\
\mathcal{F}_{r}+\mathcal{F}_{r+1}+\dots+\mathcal{F}_{r+9} & =11\mathcal{F}_{r+6}.
\end{align*}

\subsection{Fibonacci sequence with general initial conditions}

Now, we consider a Fibonacci sequence $(x_n)_{n\in \mathbb{N}}$ with general
initial conditions. We are given two complex numbers $x_0,x_1$ and we assume that
$x_{n+2}=x_{n+1}+x_n$ for any $n\in\mathbb{N}$.

Put $S_0=0$ and $S_n=\sum_{k=0}^{n-1} x_k$ for any integer $n\ge 1$. This sum
can be easily computed as before:
\begin{align*}
S_n=x_{n+1}-x_1.
\end{align*}
We propose two approaches to tackle the problem of finding a proportionality
identity of the form $S_n=Ax_m$ by appealing to simple arguments of linear algebra.
We know that the set of Fibonacci sequences (with general initial
conditions) is a vectorial space $\mathcal{E}$ of dimension~$2$.
We use decompositions of the sequence $(x_n)_{n\in \mathbb{N}}$ over two appropriate
bases of $\mathcal{E}$. Our aim is not to find all $m$ but only to derive
a convenient integer $m$ between $1$ and $n$.

\subsubsection*{\large \sl First approach}

We consistently put $\mathcal{F}_{-1}=1$, so that the sequence
$(\mathcal{F}_{n-1})_{n\in \mathbb{N}}$ is a Fibonacci sequence with initial
conditions $\mathcal{F}_{-1}=1$ and $\mathcal{F}_0=0$.
It is clear that the family $\{(\mathcal{F}_{n-1})_{n\in \mathbb{N}},
(\mathcal{F}_n)_{n\in \mathbb{N}}\}$
is a basis of $\mathcal{E}$ and that $(x_n)_{n\in \mathbb{N}}$ can be decomposed as
\begin{align*}
x_n=x_0\mathcal{F}_{n-1}+x_1\mathcal{F}_n,\quad n\in\mathbb{N}.
\end{align*}
Therefore, concerning the sum $S_n$,
\begin{align*}
S_n=x_{n+1}-x_1=x_0\mathcal{F}_n+x_1(\mathcal{F}_{n+1}-1).
\end{align*}
In particular, referring to~(\ref{rel-fibolucas}), for any integer $m\ge 1$,
\begin{align*}
S_{2m-2}=x_0\mathcal{F}_{2m-2}+x_1(\mathcal{F}_{2m-1}-1)
=\mathcal{L}_{m-1}(x_0\mathcal{F}_{m-1}+x_1\mathcal{F}_m)
=\mathcal{L}_{m-1}x_m.
\end{align*}
This approach directly supplies an answer to our problem:
for any positive integer $n$ which is a multiple of $4$ plus $2$, we have found
an integer $m=n/2+1$ and a proportionality factor $\mathcal{A}=\mathcal{L}_{n/2}$
such that $S_n=\mathcal{A}x_m$.
%
\begin{theorem}\label{theorem-fibo}
Let $(x_n)_{n\in\mathbb{N}}$ be a sequence of complex numbers satisfying the
recurrence relation $x_{n+2}=x_{n+1}+x_n$ for any $n\in\mathbb{N}$.
For any positive integer $n$ which is a multiple of\/ $4$ plus $2$, the
following identity holds:
\begin{align}\label{relS-fiboCI}
S_n=\mathcal{L}_{n/2}x_{n/2+1}.
\end{align}
\end{theorem}
%
Identity~(\ref{relS-fiboCI}) is quite analogous to~(\ref{relS-fibo}).
Moreover, it is worth noting that the factor $\mathcal{A}=\mathcal{L}_{n/2}$
is a positive integer which does not depend on the initial conditions,
according to the goal of our investigation. This approach essentially hinges
on the particular form of the recurrence relation leading to a telescopic sum,
which is quite restrictive. In order to extend Theorem~\ref{theorem-fibo}
to more general sequences, we propose a
more constructive approach which will works in a wider context.

\subsubsection*{\large \sl Second approach}

We use another basis of $\mathcal{E}$:
$\{(\varphi_1^n)_{n\in \mathbb{N}},(\varphi_2^n)_{n\in \mathbb{N}}\}$.
Then $(x_n)_{n\in \mathbb{N}}$ can be decomposed as
\begin{align*}
x_n=\alpha_1 \varphi_1^n+\alpha_2 \varphi_2^n,\quad n\in\mathbb{N},
\end{align*}
where $\alpha_1$ and $\alpha_2$ are two parameters
characterized by the equalities $\alpha_1+\alpha_2=x_0$ and $\alpha_1\varphi_1+\alpha_2\varphi_2=x_1$.
Actually, it is not necessary to determine explicitly $\alpha_1$ and $\alpha_2$.
With these settings at hand, the sum $S_n$ can be evaluated as follows:
\begin{align*}
S_n =\alpha_1\sum_{k=0}^{n-1} \varphi_1^k+\alpha_2\sum_{k=0}^{n-1} \varphi_2^k
=\alpha_1\frac{\varphi_1^n-1}{\varphi_1-1}+\alpha_2\frac{\varphi_2^n-1}{\varphi_2-1}.
\end{align*}
We would like to relate $S_n$ to one term $x_m$ of the sequence through
a simple proportionality identity independently on $\alpha_1$ and $\alpha_2$.
More precisely, is it possible to find, for any integer $n$, a factor $\mathcal{A}$
and an integer $m$ such that $S_n=\mathcal{A} x_m$ for any choice of $x_0,x_1$?

This is equivalent to finding, for any $n$, numbers $\mathcal{A}$ and $m$ such that
\begin{align*}
\forall \alpha_1,\alpha_2\in\mathbb{C},\;
\alpha_1\frac{\varphi_1^n-1}{\varphi_1-1}+\alpha_2\frac{\varphi_2^n-1}{\varphi_2-1}
=\alpha_1 \mathcal{A}\varphi_2^m+\alpha_1 \mathcal{A}\varphi_2^m
\end{align*}
or
\begin{align*}
\frac{\varphi_1^n-1}{\varphi_1-1}=\mathcal{A}\varphi_1^m \;\text{ and }\;
\frac{\varphi_2^n-1}{\varphi_2-1}=\mathcal{A}\varphi_2^m.
\end{align*}
These last equalities can be rewritten as
\begin{align*}
\varphi_1^n-1=\mathcal{A}(\varphi_1^{m+1}-\varphi_1^m) \;\text{ and }\;
\varphi_2^n-1=\mathcal{A}(\varphi_2^{m+1}-\varphi_2^m)
\end{align*}
or, since $\varphi_1^2=\varphi_1+1$ and $\varphi_2^2=\varphi_2+1$
which entail that
$\varphi_1^{m+1}-\varphi_1^m=\varphi_1^{m-1}$
and $\varphi_2^{m+1}-\varphi_2^m=\varphi_2^{m-1}$,
\begin{align}\label{proportionA-inter}
\varphi_1^n-1=\mathcal{A}\varphi_1^{m-1} \;\text{ and }\;
\varphi_2^n-1=\mathcal{A}\varphi_2^{m-1}.
\end{align}
Hence, we are searching for $\mathcal{A}$ and $m$ such that
\begin{align}\label{proportionA}
\mathcal{A}=\frac{\varphi_1^n-1}{\varphi_1^{m-1}}
=\frac{\varphi_2^n-1}{\varphi_2^{m-1}}.
\end{align}
Consequently, $m$ solves the equation
\begin{align*}
\left(\frac{\varphi_1}{\varphi_2}\right)^{\!m-1}
=\frac{\varphi_1^n-1}{\varphi_2^n-1}.
\end{align*}
Substituting $\varphi_2=-1/\varphi_1$ in the foregoing equation
yields
\begin{align*}
\varphi_1^{2m-n-2}=(-1)^m\frac{\varphi_1^n-1}{\varphi_1^n-(-1)^n}.
\end{align*}
When $n$ is odd, it is difficult to extract $m$ and it seems hopeless
to obtain an integer. When $n$ is even instead, this is equivalent to
\begin{align*}
\varphi_1^{2m-n-2}=(-1)^m.
\end{align*}
Since $\varphi_1>0$, the integer $m$ must be even too and $n=2m-2$.
As a result, $n$ must be a multiple of $4$ plus $2$ and
$m=n/2+1$. Next, by~(\ref{proportionA}), the coefficient $\mathcal{A}$ writes
\begin{align*}
\mathcal{A}=\frac{\varphi_1^n-1}{\varphi_1^{n/2}}=\varphi_1^{n/2}-\varphi_1^{-n/2}
=\varphi_1^{n/2}+\varphi_2^{n/2}=\mathcal{L}_{n/2}.
\end{align*}
We just retrieved the solution~(\ref{relS-fiboCI}) obtained by the first approach.
%
\begin{remark}
Let us remind that if an equality between two ratios
$\displaystyle f\!\!:=\frac{a_1}{b_1}=\frac{a_2}{b_2}$ holds, then
$\displaystyle f=\frac{a_1-a_2}{b_1-b_2}$; this is quite an elementary fact!
This yields here
\begin{align*}
\mathcal{A}=\frac{\varphi_1^n-\varphi_2^n}{\varphi_1^{m-1}-\varphi_2^{m-1}}
=\frac{\mathcal{F}_n}{\mathcal{F}_{m-1}}
\end{align*}
and we retrieve for $m=n/2+1$,
$\mathcal{A}=\mathcal{F}_n/\mathcal{F}_{n/2}=\mathcal{L}_{n/2}$.
\end{remark}
%

The second approach allows us to tackle our problem for more
general sequences. We employ it in the next section.

A quite similar analysis can be carried out
\textit{mutatis mutandis} for the sequence $(x_n)_{n\in\mathbb{N}}$ defined by
general initial conditions $x_0,x_1\in\mathbb{C}$ and the recurrence relation
$x_{n+2}=-x_{n+1}+x_n$ for any $n\in\mathbb{N}$.
The numbers $\varphi_1$ and $\varphi_2$ have to be changed into
the roots $\psi_1$ and $\psi_2$ of the equation $\psi^2=-\psi+1$.
Explicitly, $\psi_1=-\varphi_1$ and $\psi_2=-\varphi_2$.
The sequence corresponding to the initial values $0$ and $1$ consists of
the numbers $(\psi_1^n-\psi_2^n)/(\psi_1-\psi_2)$ which are nothing
but $(-1)^{n-1}\mathcal{F}_n$. Since $\psi_1^{m+1}-\psi_1^{m}=-\psi_1^{m+2}$ and
$\psi_2^{m+1}-\psi_2^{m}=-\psi_2^{m+2}$, the analogous of~(\ref{proportionA-inter}) is
\begin{align*}
\psi_1^n-1=-\mathcal{A}'\psi_1^{m+2} \;\text{ and }\;
\psi_2^n-1=-\mathcal{A}'\psi_2^{m+2}
\end{align*}
from which we extract
\begin{align*}
\left(\frac{\psi_1}{\psi_2}\right)^{\!m+2}=\frac{\psi_1^n-1}{\psi_2^n-1}
\end{align*}
or, equivalently when $n$ is even, $\varphi_1^{2m-n+4}=(-1)^{m-1}$.
We see that the integer $m$ must be odd and $n=2m+4$. So,
the integer $n$ must be a multiple of $4$ plus $2$ and
$m=n/2-2$. Next,
\begin{align*}
\mathcal{A}'=-\frac{\psi_1^n-1}{\psi_1^{n/2}}
=\psi_1^{-n/2}-\psi_1^{n/2}=\varphi_1^{n/2}+\varphi_2^{n/2}
=\mathcal{L}_{n/2}.
\end{align*}
We observe that $\mathcal{A}'=\mathcal{A}$. As a result, we can state the result below.
%
\begin{theorem}
Let $(x_n)_{n\in\mathbb{N}}$ be a sequence satisfying the recurrence relation
$x_{n+2}=-x_{n+1}+x_n$ for any $n\in\mathbb{N}$.
For any positive integer $n$ which is a multiple of\/ $4$ plus $2$, the
following identity holds:
\begin{align*}
S_n=\mathcal{L}_{n/2}x_{n/2-2}.
\end{align*}
\end{theorem}
%

As an example of a trick, for $n=10$, you choose two numbers $x_0$
and $x_1$ and you successively write
\begin{eqnarray*}
&x_2=x_0-x_1,\; x_3=-x_0+2x_1,\; x_4=2x_0-3x_1,\; x_5=-3x_0+5x_1,&
\\
&x_6=5x_0-8x_1,\; x_7=-8x_0+13x_1,\; x_8=13x_0-21x_1,\; x_9=-21x_0+34x_1.&
\end{eqnarray*}
Finally, you compute the sum of all the numbers from $x_0$ to $x_9$:
$S=x_0+x_1+\dots+x_9$. You find that
\begin{align*}
S=-11x_0+22x_1=11x_3.
\end{align*}
We guess that this trick is simple enough to be exploited in magic...

\section{The second order linear recurrence}\label{section-Lucas}

Let us fix two complex numbers $a,b$ and let us consider the sequence
$(x_n)_{n\in \mathbb{N}}$ defined by two initial conditions
$x_0,x_1\in\mathbb{C}$ and the recurrence relation $x_{n+2}=ax_{n+1}+bx_n$ for
any $n\in\mathbb{N}$. When the parameters $a$ and $b$ are integers,
$(x_n)_{n\in \mathbb{N}}$ is called a Lucas sequence; see~\cite{lucasseq}
and~\cite{wikilucasseq}.

Let us introduce $r_1$ and $r_2$ the complex roots of the quadratic equation
$r^2=ar+b$. We assume that $a^2+4b\neq 0$ in order to have two distinct roots
($r_1\neq r_2$), that $b\neq 0$ and $a+b\neq 1$ in order that $0$ and $1$
are not roots ($r_1,r_2\notin\{0,1\}$). The numbers $r_1$ and $r_2$ satisfy
$r_1r_2=-b$.

The explicit expression of $x_n$ is then given by
$x_n=\alpha_1 r_1^n+\alpha_2 r_2^n$ where $\alpha_1$ and $\alpha_2$ are
two parameters such that $\alpha_1+\alpha_2=x_0$ and
$\alpha_1r_1+\alpha_2r_2=x_1$. It will not be necessary to determine
$\alpha_1$ and $\alpha_2$.
Put $S_0=0$ and $S_n=\sum_{k=0}^{n-1} x_k$ for $n\ge 1$. We have
\begin{align}\label{sum}
S_n=\alpha_1\frac{r_1^n-1}{r_1-1}+\alpha_2\frac{r_2^n-1}{r_2-1}.
\end{align}
As previously, is it possible to find, for any integer $n$, a coefficient $A$
and an integer $m$ such that $S_n=A x_m$ for any choice of $x_0$ and $x_1$?

\subsection{Solving the problem}

The problem is equivalent to finding, for any $n$, numbers $A$ and $m$ such that
\begin{align*}
r_1^n-1=A(r_1-1)r_1^m \;\text{ and }\; r_2^n-1=A(r_2-1)r_2^m.
\end{align*}
So, we get
\begin{align}\label{proportionA2}
A=\frac{r_1^n-1}{(r_1-1)r_1^m}=\frac{r_2^n-1}{(r_2-1)r_2^m}
\end{align}
which entails
\begin{align*}
\left(\frac{r_1}{r_2}\right)^{\!m}=\frac{r_2-1}{r_1-1}\,\frac{r_1^n-1}{r_2^n-1}.
\end{align*}
Substituting $r_2=-b/r_1$ in the foregoing equation yields
\begin{align}\label{equation-difficult}
r_1^{2m-n+1}=(-b)^m\frac{r_1+b}{r_1-1}\,\frac{r_1^n-1}{r_1^n-(-b)^n}.
\end{align}
When $b\neq \pm1$, Equation~(\ref{equation-difficult}) is difficult to solve
and it seems hopeless to find an integer $m$ satisfying it; we shall not
go further in this direction. We now focus precisely on the cases
$b=-1$ and $b=1$.

\subsubsection*{\large \sl The case $b=-1$}

Suppose that $b=-1$. The assumptions $a^2+4b\neq 0$ and $a+b\neq 0$
force $a$ to be chosen different from $2$ and $-2$.
Equation~(\ref{equation-difficult}) reads $r_1^{2m-n+1}=1$ which implies that
$n=2m+1$. Hence, $n$ must be odd and $m=(n-1)/2$.

Next, we have to determine the coefficient $A$.
Using the proportionality rule stipulating that if
$\alpha:=\alpha_1=\alpha_2$, then $\alpha=[(r_1-1)\alpha_1
-(r_2-1)\alpha_2]/(r_1-r_2)$, and the fact that $r_1r_2=-b=1$,
we extract from~(\ref{proportionA2}) the alternative form
\begin{align}\label{proportionA-interbis}
A=\frac{1}{r_1-r_2}\left(\frac{r_1^n-1}{r_1^m}-\frac{r_2^n-1}{r_2^m}\right)
=\frac{r_1^{n-m}-r_2^{n-m}}{r_1-r_2}+\frac{r_1^m-r_2^m}{r_1-r_2}.
\end{align}
By introducing the particular sequence $(u_n)_{n\in \mathbb{N}}$
satisfying the same recurrence relation as $(x_n)_{n\in \mathbb{N}}$ with
initial conditions $u_0=0$ and $u_1=1$, whose expression is
$\displaystyle u_n=(r_1^n-r_2^n)/(r_1-r_2),$
we see that the number $A$ can be written as $A=u_{n-m}+u_m=u_{m+1}+u_m$, that is
\begin{align}\label{proportionA2bis}
A=u_{(n+1)/2}+u_{(n-1)/2}.
\end{align}
In Section~\ref{subsection-solving-cont}, we propose a closed expression
for $A$ by means of $a$.

Let us have a look at the cases $a=\pm2$ which were excluded from the previous
analysis.
\begin{itemize}
\item
When $a=2$, $x_n$ is of the form $x_n=\alpha_1 n+\alpha_2$ for certain parameters
$\alpha_1$ and $\alpha_2$, $u_n=n$ for any $n\in\mathbb{N}$ and, if $n$ is odd,
\begin{align*}
S_n=\alpha_1\,\frac{n(n-1)}{2}+n\alpha_2
=n\!\left(\alpha_1\,\frac{n-1}{2}+\alpha_2\right)=nx_{(n-1)/2}.
\end{align*}
We find that the equality $S_n=Ax_{(n-1)/2}$ holds with $A=n$
which coincides with $u_{(n+1)/2}+u_{(n-1)/2}$.
\item
When $a=-2$, we similarly have
$x_n=(-1)^n(\alpha_1 n+\alpha_2)$, $u_n=(-1)^{n+1}n$ for any $n\in\mathbb{N}$ and,
if $n=2m+1$,
\begin{align*}
S_n=S_{2m+1}=\alpha_1 m+\alpha_2=(-1)^m x_m=(-1)^{(n-1)/2} x_{(n-1)/2}.
\end{align*}
Then $S_n=Ax_{(n-1)/2}$ with $A=(-1)^{(n-1)/2}$. This $A$ coincides with
$u_{(n+1)/2}+u_{(n-1)/2}$ and (\ref{proportionA2bis}) holds in this case also.
\end{itemize}
We conclude as follows.
%
\begin{theorem}\label{theorem-gene}
Fix a complex number $a$. Let $(x_n)_{n\in\mathbb{N}}$ be a sequence satisfying
the recurrence relation $x_{n+2}=ax_{n+1}-x_n$ for any $n\in\mathbb{N}$
and $(u_n)_{n\in\mathbb{N}}$ be the sequence satisfying the same relation such
that $u_0=0$ and $u_1=0$.
For any positive and odd integer $n$, the following proportionality identity holds:
\begin{align*}
S_n=(u_{(n+1)/2}+u_{(n-1)/2})x_{(n-1)/2}.
\end{align*}
\end{theorem}
%
\begin{example}
Let us write the first numbers of the sequence:
\begin{align*}
x_2&=-x_0+ax_1,
\\
x_3&=-ax_0+(a^2-1)x_1,
\\
x_4&=(-a^2+1)x_0+(a^3-2a)x_1,
\\
x_5&=(-a^3+2a)x_0+(a^4-3a^2+1)x_1,
\\
x_6&=(-a^4+3a^2-1)x_0+(a^5-4a^3+3a)x_1,
\\
x_7&=(-a^5+4a^3-3a)x_0+(a^6-5a^4+6a^2-1)x_1,\;
\\
x_8&=(-a^6+5a^4-6a^2+1)x_0+(a^7-6a^5+10a^3-4a)x_1,
\\
x_9&=(-a^7+6a^5-10a^3+4a)x_0+(a^8-7a^6+15a^4-10a^2+1)x_1,
\\
x_{10}&=(-a^8+7a^6-15a^4+10a^2-1)x_0+(a^9-8a^7+21a^5-20a^3+5a)x_1.
\end{align*}
When choosing $x_0=0$ and $x_1=1$, we obtain $u_5=a^4-3a^2+1$
as well as $u_6=a^5-4a^3+3a$. The sum of all numbers from $x_0$ to $x_{10}$
(\/$n=11$) is given by
\begin{align*}
S_{11}&=(-a^8-a^7+6a^6+5a^5-11a^4-7a^3+6a^2+2a)x_0
\\
&\phantom{=}\;+(a^9+a^8-7a^7-6a^6+16a^5+11a^4-13a^3-6a^2+3a+1)x_1
\end{align*}
The reader can easily check that
\begin{align*}
S_{11}=(a^5+a^4-4a^3-3a^2+3a+1)[(-a^3+2a)x_0+(a^4-3a^2+1)x_1]=(u_5+u_6)x_5.
\end{align*}
\end{example}
%

\subsubsection*{\large \sl The case $b=1$}

Let us now consider the case where $b=1$. We are dealing with
the recurrence relation $x_{n+2}=ax_{n+1}+x_n$. Assume that $m$ and $n$ are even.
Equation~(\ref{equation-difficult}) reads
\begin{align}\label{equation-simple}
r_1^{2m-n+1}=\frac{r_1+1}{r_1-1}.
\end{align}
Fix a positive integer $p$ and choose a complex root $\rho_1$ of
the equation $r^{p+1}-r^p-r-1=0$. Set $\rho_2=-1/\rho_1$.
If we choose $a=\rho_1+\rho_2$, the roots $\rho_1$ and $\rho_2$
satisfy $\rho_1^2=a\rho_1+1$ and $\rho_2^2=a\rho_2+1$.
In order to use our analysis, the equation $r^2=ar+1$ must have two
different roots, that is,
we must choose $\rho_1$ different from $\pm i$. The numbers $\pm i$
are roots when the power $p$ is a multiple of $4$ minus $1$.
Since $\rho_1^p=(\rho_1+1)/(\rho_1-1)$, the number $\rho_1$
solves Equation~(\ref{equation-simple}) if and only if $\rho_1^{2m-n+1}=\rho_1^p$
which is clearly satisfied when $n=2m+1-p$, or $m=(n+p-1)/2$.
Taking into account the assumptions that $m$ and $n$ are even,
$p$ must be chosen \textit{a priori} odd such that $(n+p-1)$ be a multiple
of $4$.

Finally, by introducing the particular sequence $(u_n)_{n\in \mathbb{N}}$
satisfying the same recurrence relation as $(x_n)_{n\in \mathbb{N}}$ with
initial conditions $u_0=0$ and $u_1=1$, it is quite easy to see that the
expression $A=u_{n-m}+u_m$ for the proportionality coefficient provided
in the case $b=-1$ still holds, due to the fact that $m$ is even
and $\rho_1\rho_2=-1$. Therefore,
\begin{align}\label{proportionA2ter}
A=u_{(n-p+1)/2}+u_{(n+p-1)/2}.
\end{align}
%
\begin{remark}
The equation $r^{p+1}-r^p-r-1=0$ admits at least one positive real root:
indeed, it is equivalent to $r^p=\frac{r+1}{r-1}$, the roots of which are
the intersection of the curves of the maps $x\mapsto x^p$ and
$x\mapsto\frac{x+1}{x-1}$. More precisely, an elementary study shows that
the equation has exactly one positive root (which lies in the interval $(1,2)$)
when $p$ is even and two roots (one is positive--lying in $(1,2.5)$--and the
other is the opposite of the inverse of the first one) if $p$ is odd.
\end{remark}
%
\begin{remark}
Referring to the first equality of~(\ref{proportionA-interbis}) and using
the facts that $m$ and $n$ are even and that $\rho_1\rho_2=-1$,
we can derive another expression of $A$:
\begin{align*}
A & =\frac{1}{\rho_1-\rho_2}\left(\frac{\rho_1^n-1}{\rho_1^m}-
\frac{\rho_2^n-1}{\rho_2^m}\right)
\\
& =\rho_1^{n/2-m}\,\frac{\rho_1^{n/2}-\rho_1^{-n/2}}{\rho_1-\rho_2}
-\rho_2^{n/2-m}\,\frac{\rho_2^{n/2}-\rho_2^{-n/2}}{\rho_1-\rho_2}
\\
& =\rho_1^{n/2-m}\,\frac{\rho_1^{n/2}-(-1)^{n/2}\rho_2^{n/2}}{\rho_1-\rho_2}
-\rho_2^{n/2-m}\,\frac{\rho_2^{n/2}-(-1)^{n/2}\rho_1^{n/2}}{\rho_1-\rho_2}
\\
& =\frac{\left(\rho_1^{n/2-m}+(-1)^{n/2}\rho_2^{n/2-m}\right)\!\!
\left(\rho_1^{n/2}-(-1)^{n/2}\rho_2^{n/2}\right)}{\rho_1-\rho_2}
\\
& =\frac{\left(\rho_1^{m-n/2}+(-1)^{n/2}\rho_2^{m-n/2}\right)\!\!
\left(\rho_1^{n/2}-(-1)^{n/2}\rho_2^{n/2}\right)}{\rho_1-\rho_2}
\\
& =\frac{\left(\rho_1^{(p-1)/2}+(-1)^{n/2}\rho_2^{(p-1)/2}\right)\!\!
\left(\rho_1^{n/2}-(-1)^{n/2}\rho_2^{n/2}\right)}{\rho_1-\rho_2}.
\end{align*}
As a byproduct, by setting $u_q=(\rho_1^q-\rho_2^q)/(\rho_1-\rho_2)$
and $v_q=\rho_1^q+\rho_2^q$ for any $q\in\mathbb{N}$,
\begin{align}\label{proportionA2quater}
A=\begin{cases}
v_{(p-1)/2}\,u_{n/2} & \text{if $n$ is a multiple of $4$,}
\\[.2ex]
u_{(p-1)/2}\,v_{n/2} & \text{if $n$ is a multiple of $4$ plus $2$.}
\end{cases}
\end{align}
\end{remark}
%
We conclude the case $b=1$ as follows.
\begin{theorem}\label{theorem-genebis}
Fix a positive and odd integer $p$ and let $\rho$ be a root of the equation
$r^{p+1}-r^p-r-1=0$ different from $\pm i$. Set $a=\rho-1/\rho$ and
let us introduce a sequence $(x_n)_{n\in\mathbb{N}}$ satisfying
the recurrence relation $x_{n+2}=ax_{n+1}+x_n$ for any $n\in\mathbb{N}$
and $(u_n)_{n\in\mathbb{N}}$ and $(v_n)_{n\in\mathbb{N}}$ the sequences
satisfying the same relation such that $u_0=0$ and $u_1=1$, $v_0=2$ and
$v_1=a$ respectively.
For any positive and even integer $n$ such that $(n+p-1)$ is a multiple
of\/ $4$, the following identity holds:
\begin{align*}
S_n & =(u_{(n-p+1)/2}+u_{(n+p-1)/2})\,x_{(n+p-1)/2}
=\begin{cases}
v_{(p-1)/2}\,u_{n/2}\,x_{(n+p-1)/2} & \text{if $n$ is a multiple of $4$,}
\\[.2ex]
u_{(p-1)/2}\,v_{n/2}\,x_{(n+p-1)/2} & \text{if $n$ is a multiple of $4$ plus $2$.}
\end{cases}
\end{align*}
\end{theorem}
%
\begin{examples}\label{examples}
\begin{itemize}
\item
For $p=1$, the corresponding equation reads $r^2=2r+1$ whose solutions are
$\big(1+\sqrt2\,\big)$ and $\big(1-\sqrt2\,\big)$. Hence $a=2$, the recurrence relation
for the sequence $(x_n)_{n\in\mathbb{N}}$ writes $x_{n+2}=2x_{n+1}+x_n$ and
the proportionality identity reads, for any integer $n$ which is a multiple of\/ $4$,
\begin{align*}
S_n=2\,u_{n/2}\,x_{n/2}.
\end{align*}
As an example of a possible trick, for $n=12$, choose two numbers $x_0$ and $x_1$
and compute the $10$ next numbers following the recurrence relation
$x_{n+2}=2x_{n+1}+x_n$:
\begin{eqnarray*}
&x_2=x_0+2x_1,\; x_3=2x_0+5x_1,\; x_4=5x_0+12x_1,\; x_5=12x_0+29x_1,&
\\
&x_6=29x_0+70x_1,\; x_7=70x_0+169x_1,\; x_8=169x_0+408x_1,&
\\
&x_9=408x_0+985x_1,\;x_{10}=985x_0+2378x_1,\; x_{11}=2378x_0+5741x_1.&
\end{eqnarray*}
The strategic number $u_6$ is obtained from $x_6$ by substituting therein the values
$0$ and $1$ into $x_0$ and $x_1$: $u_6=70$.
Evaluate now the sum of all these numbers $S=x_0+x_1+\dots+x_{11}$. You find that
\begin{align*}
S=4060 x_0+ 9800 x_1=140(29 x_0+70 x_1)=140 x_6=2 u_6 x_6.
\end{align*}

\item
For $p=3$, the corresponding equation writes $r^4-r^3-r-1=0$ whose solutions are
$i$, $-i$, $\varphi_1$ and $\varphi_2$. The roots $\varphi_1$ and $\varphi_2$
lead to $a=1$ and we retrieve the Fibonacci sequence. The roots $i$
and $-i$ are excluded from our analysis. The sequences
$(u_n)_{n\in\mathbb{N}}$ and $(v_n)_{n\in\mathbb{N}}$ are the Fibonacci and Lucas sequences
$(\mathcal{F}_n)_{n\in\mathbb{N}}$ and $(\mathcal{L}_n)_{n\in\mathbb{N}}$ respectively,
and, for any integer $n$ which is a multiple of\/ $4$ plus\/ $2$,
\begin{align*}
S_n=(\mathcal{F}_{n/2+1}+\mathcal{F}_{n/2-1})\,x_{n/2+1}=\mathcal{L}_{n/2}\,x_{n/2+1}.
\end{align*}

\item
For $p=5$, the corresponding equation writes $r^6-r^5-r-1=0$ which admits two
real solutions. Maple supplies the following estimates: $-0.704$ and $1.420$
(we have $-1/1.420\approx -0.704$).
These both give the same possible value of $a$: $a\approx 0.715$;
the recurrence relation for the corresponding sequence $(x_n)_{n\in\mathbb{N}}$
reads $x_{n+2}=ax_{n+1}+x_n$, and the proportionality identity reads,
for any $n$ which is a multiple of\/ $4$,
\begin{align*}
S_n=(u_{n/2+2}+u_{n/2-2})x_{n/2+2}=(a^2+2)u_{n/2}x_{n/2+2}.
\end{align*}
We postpone to the appendix additional information about this example.
In particular, we check this formula for $n=8$.
\end{itemize}
\end{examples}
%
In the table below, we collect estimates for certain coefficients $a$
related to small values of the odd integer $p$ supplied by Maple.
The number $\rho$ is the positive root of the polynomial
$r^{p+1}-r^p-r-1$, $a=\rho-1/\rho$ determines the sequence for
which a proportionality identity holds,
and $A$ is the proportionality coefficient.
\newcommand{\bloc}[1]{\raisebox{-.3ex}{$#1$}}
$$
\begin{array}{|@{\hspace{.8em}}c@{\hspace{.8em}}||@{\hspace{.8em}}
c@{\hspace{.8em}}|c@{\hspace{.6em}}||c@{\hspace{.5em}}|c@{\hspace{.8em}}|}
\hline
p  & \rho  & a     & m     & \bloc{A}
\\
\hline\hline
\bloc{1}  & \bloc{2.414} & \bloc{2}     & \bloc{n/2} & \bloc{2\,u_{n/2}}
\\
\hline
\bloc{3}  & \bloc{1.618} & \bloc{1}     & \bloc{n/2+1} & 
v_{n/2}
\\
\hline
\bloc{5}  & \bloc{1.420} & \bloc{0.715} & \bloc{n/2+2} & 
\bloc{(a^2+2)u_{n/2}}
\\
\hline
\bloc{7}  & \bloc{1.325} & \bloc{0.570} & \bloc{n/2+3} & 
\bloc{(a^2+1)v_{n/2}}
\\
\hline
\bloc{9}  & \bloc{1.268} & \bloc{0.479} & \bloc{n/2+4} & 
\bloc{(a^4+4a^2+2)u_{n/2}}
\\
\hline
\bloc{11} & \bloc{1.230} & \bloc{0.416} & \bloc{n/2+5} & 
\bloc{(a^4+3a^2+1)v_{n/2}}
\\
\hline
\bloc{13} & \bloc{1.202} & \bloc{0.370} & \bloc{n/2+6} & 
\bloc{(a^6+6a^4+9a^2+2)u_{n/2}}
\\
\hline
\bloc{15} & \bloc{1.181} & \bloc{0.334} & \bloc{n/2+7} & 
\bloc{(a^6+5a^4+6a^2+1)v_{n/2}}
\\
\hline
\bloc{17} & \bloc{1.164} & \bloc{0.305} & \bloc{n/2+8} & 
\bloc{(a^8+8a^6+20a^4+16a^2+2)u_{n/2}}
\\
\hline
\bloc{19} & \bloc{1.150} & \bloc{0.281} & \bloc{n/2+9} & 
\bloc{(a^8+7a^6+15a^4+10a^2+1)v_{n/2}}
\\
\hline
\end{array}
$$

\subsection{Solving the problem (continued)}\label{subsection-solving-cont}

The foregoing analysis suggests considering only the family of sequences
$(x_n)_{n\in \mathbb{N}}$ following one of the recurrence relations
$x_{n+2}=ax_{n+1}-x_n$ or $x_{n+2}=ax_{n+1}+x_n$ for any $n\in\mathbb{N}$,
where $a\in\mathbb{C}$ is a fixed parameter. Actually, we essentially restrict
ourselves to the first relation since the second one only concerns certain
specific values of $a$. We briefly examine the second one at the end of this section.

For the sequences satisfying $x_{n+2}=ax_{n+1}-x_n$,
we evaluate in this part the coefficient $A$ given
by~(\ref{proportionA2bis}) by means of the parameter $a$. We guess that
the reader has certainly recognized the recurrence relation of the famous
Chebyshev polynomials. We make clear this connection in Section~\ref{section-cheb}.
We propose here two representations for $A$.

As usual, $(u_n)_{n\in \mathbb{N}}$ denotes the sequence satisfying the
same relation as $(x_n)_{n\in \mathbb{N}}$ with initial conditions $u_0=0$
and $u_1=1$.
We recall the expression of $u_n$: $u_n=(r_1^n-r_2^n)/(r_1-r_2)$ where $r_1$ and $r_2$ are the
roots of $r^2=ar-1$, as well as that of $A$: $A=u_{(n+1)/2}+u_{(n-1)/2}$.

\subsubsection*{\large \sl First approach}

A first way of computing $A$ consists in using generating functions.
The generating function of $(u_n)_{n\in \mathbb{N}}$
is defined, for any complex number $z$ such that $|z|<\min(|r_1|,|r_2|)$, by
$G(z)=\sum_{n=0}^{\infty} u_nz^n$. It satisfies the following equation:
\begin{align*}
G(z) & =u_0+u_1z+\sum_{n=2}^{\infty} (au_{n-1}-u_{n-2})z^n
=z+az\sum_{n=2}^{\infty} u_{n-1}z^{n-1}-z^2\sum_{n=2}^{\infty}u_{n-2}z^{n-2}
\\
& =z+az\sum_{n=1}^{\infty} u_nz^n-z^2\sum_{n=0}^{\infty}u_nz^n
=z+(az-z^2)G(z).
\end{align*}
We immediately extract
\begin{align*}
G(z)=\frac{z}{1-(az-z^2)}.
\end{align*}
Now, we expand $G(z)$ into a power series by using the elementary
identity $1/(1-\zeta)=\sum_{\ell=0}^{\infty}\zeta^{\ell}$:
\begin{align*}
G(z) & =\sum_{\ell=0}^{\infty} z(az-z^2)^{\ell}
=\sum_{\ell=0}^{\infty} \sum_{k=0}^{\ell} (-1)^k\binom{\ell}{k} a^{\ell-k}z^{k+\ell+1}
=\sum_{n=1}^{\infty} \Bigg(\sum_{k=0}^{\lfloor (n-1)/2\rfloor}
(-1)^k\binom{n-k-1}{k} a^{n-2k-1}\Bigg)z^n.
\end{align*}
By identification, we derive the following representation of $u_n$:
\begin{align*}
u_n=\sum_{k=0}^{\lfloor (n-1)/2\rfloor} (-1)^k\binom{n-k-1}{k} a^{n-2k-1}.
\end{align*}
As a result, with $m=(n-1)/2$,
\begin{align}\label{proportionA3}
A=u_{m+1}\!+u_m=\!\sum_{k=0}^{\lfloor m/2\rfloor}
\!(-1)^k\!\left[\!\binom{m-k-1}{k} a^{m-2k-1}
+\binom{m-k}{k} a^{m-2k}\right]\!.
\end{align}
In Formula~(\ref{proportionA3}), we adopt the convention that
$\binom{\ell}{k}=0$ if $k>\ell$.
We observe that the coefficient $A$ is a polynomial of the parameter $a$
of the sequence $(x_n)_{n\in\mathbb{N}}$. This fact is actually not at all
surprising and we detail this in Section~\ref{subsection-cheb}.

Very similar computations would lead, for the sequence $(u_n)_{n\in \mathbb{N}}$
satisfying $u_0=0$, $u_1=1$ and the recurrence relation $u_{n+2}=au_{n+1}+u_n$,
to
\begin{align*}
u_n=\sum_{k=0}^{\lfloor (n-1)/2\rfloor} \binom{n-k-1}{k} a^{n-2k-1}.
\end{align*}

\subsubsection*{\large \sl Second approach}

A second expression for $A$ can be found. For this, we examine several cases
according to the values of $a$ when $a\in\mathbb{R}\backslash\{-2,2\}$ by using
the explicit expressions of the roots of $r^2=ar-1$.
\begin{itemize}
\item
If $|a|<2$, we can write $a=2\cos\theta$ for a certain
$\theta\in\,(0,\pi)$ and the roots of $r^2=ar-1$ are
$e^{i\theta}$ and $e^{-i\theta}$.
Then $u_n=\sin(n\theta)/\sin\theta$ and
\begin{align*}
A=\frac{\sin((n+1)\theta/2)+\sin((n-1)\theta/2)}{\sin\theta}
=\frac{\sin(n\theta/2)}{\sin(\theta/2)}.
\end{align*}

\item
If $a>2$, we can write $a=2\cosh\theta$ for a certain $\theta>0$ and
the roots of $r^2=ar-1$ are $e^{\theta}$ and $e^{-\theta}$.
Then $u_n=\sinh(n\theta)/\sinh\theta$ and
\begin{align*}
A=\frac{\sinh(n\theta/2)}{\sinh(\theta/2)}.
\end{align*}

\item
If $a<-2$, we can write $a=-2\cosh\theta$ for a certain $\theta>0$ and
the roots of $r^2=ar-1$ are $-e^{\theta}$ and $-e^{-\theta}$. Then
$u_n=(-1)^{n-1}\sinh(n\theta)/\sinh\theta$ and
\begin{align*}
A=(-1)^{(n-1)/2}\frac{\cosh(n\theta/2)}{\cosh(\theta/2)}.
\end{align*}
\end{itemize}

%
\begin{examples}

\begin{itemize}
\item
For $a=1$, the corresponding value of the angle $\theta$ is $\pi/3$. In this case,
we have successively $x_2=x_1-x_0$, $x_3=-x_0$, $x_4=-x_1$, $x_5=x_0-x_1=-x_2$,
$x_6=x_0$. We observe that the sequence $(x_n)_{n\in\mathbb{N}}$ is periodic:
$x_{t+6}=x_t$ for any $t\in\mathbb{N}$. Next, the successive values of $S_n$ are
$S_0=0$, $S_1=x_0$, $S_2=x_0+x_1$, $S_3=2x_1$, $S_4=2x_1-x_0$, $S_5=x_1-x_0$,
$S_6=0$. The sequence $(S_n)_{n\in\mathbb{N}}$ is also periodic: $S_{t+6}=S_t$.
We can see that, for any $s,t\in\mathbb{N}$, $S_{6t+1}=x_0=x_{6s}$,
$S_{6t+3}=2x_1=2x_{6s+1}$ and $S_{6t+5}=x_2=x_{6s+2}$.

\item
For $a=-1$, the corresponding value of $\theta$ is $2\pi/3$. In this case,
we have $x_2=-x_0-x_1$, $x_3=x_0$ and $S_0=0$, $S_1=x_0$, $S_2=x_0+x_1=-x_2$, $S_3=0$.
Thus, the sequences $(x_n)_{n\in\mathbb{N}}$ and $(S_n)_{n\in\mathbb{N}}$ are periodic:
$x_{t+3}=x_t$ and $S_{t+3}=S_t$ for any $t\in\mathbb{N}$.
We see that $S_{3t+1}=x_{3s}$ and $S_{3t+2}=-x_{3s+2}$.

\item
For $a=0$, the corresponding value of $\theta$ is $\pi/2$. In this case,
we have $x_2=-x_0$, $x_3=-x_1$,
$x_4=-x_2=x_0$ and $S_0=0$, $S_1=x_0$, $S_2=x_0+x_1$, $S_3=x_1$, $S_4=0$.
In this case again, the sequences are periodic and
$S_{4t+1}=x_0=x_{4s}$ and $S_{4t+3}=x_1=x_{4s+1}$ for any $s,t\in\mathbb{N}$.
\end{itemize}
\end{examples}
%
\begin{remark}
The sequence $(x_n)_{n\in\mathbb{N}}$ is periodic if and only if the angles
$\theta$ and $\pi$ are commensurable. In this case, in view of~(\ref{sum}),
$(S_n)_{n\in\mathbb{N}}$ is clearly periodic.
\end{remark}

We end up this section by examining briefly the recurrence relation
$x_{n+2}=ax_{n+1}+x_n$. We recall that we found a solution to our
problem only for certain values of $a$ (see Theorem~\ref{theorem-genebis}).
We can write $a=2\sinh\theta$ for a certain $\theta\in\mathbb{R}$ and the roots
of the equation $r^2=ar+1$ are $e^{\theta}$ and $-e^{-\theta}$.
Then $u_n=\sinh(n\theta)/\cosh\theta$ if $n$ is even and
the expression~(\ref{proportionA2ter}) of the coefficient $A$ becomes
\begin{align*}
A & =\frac{\sinh((n+p-1)\theta/2)+\sinh((n-p+1)\theta/2)}{\cosh\theta}
=2\,\frac{\sinh(n\theta/2)\cosh((p-1)\theta/2)}{\cosh\theta}
\end{align*}
where $n$ and $p$ are positive integers such that $n$ is even, $p$ is
odd and $(n+p-1)$ is a multiple of $4$.

\section{Link with Fibonacci, Lucas and Chebyshev polynomials}\label{section-cheb}

In this part, we make clear the link between the sequences defined by the recurrence
relations $x_{n+2}=ax_{n+1}\pm x_n$ that emerged in the foregoing section
and the famous Fibonacci, Lucas and Chebyshev polynomials.

\subsection{Link with Chebyshev polynomials}\label{subsection-cheb}

We recall the definition of the families of Chebyshev polynomials
$(T_n)_{n\in\mathbb{N}}$ (first kind) and $(U_n)_{n\in\mathbb{N}}$ (second kind)
through the expansion of the trigonometric quantities $\cos(n\theta)$
and $\sin(n\theta)$ (or hyperbolic quantities $\cosh(n\theta)$ and
$\sinh(n\theta)$) as polynomials of $\cos\theta$ (or $\cosh(n\theta)$):
\begin{align*}
\cos(n\theta)=T_n(\cos\theta) \;\text{ and }\;
\sin(n\theta)=\sin(\theta)U_{n-1}(\cos\theta).
\end{align*}
We refer to~\cite{chebpol1}, \cite{chebpol2} and~\cite{wikicheb}.
They are characterized by the initial conditions $T_0(x)=1$,
$T_1(x)=x$, $U_0(x)=1$, $U_1(x)=2x$ for any $x\in\mathbb{C}$ and the relationships,
for any $n\in\mathbb{N}$ and $x\in\mathbb{C}$,
\begin{align*}
T_{n+2}(x)=2xT_{n+1}(x)-T_n(x) \;\text{ and }\; U_{n+2}(x)=2xU_{n+1}(x)-U_n(x).
\end{align*}
Additionally, we put $U_{-1}(x)=0$.
Regarding the sequence we were dealing with in Section~\ref{subsection-solving-cont},
it is convenient to introduce $\tilde{T}_n(x)=2\,T_n(x/2)$ and
$\tilde{U}_n(x)=U_{n-1}(x/2)$.
In the literature, $\tilde{T}_n$ and $\tilde{U}_n$ are usually
denoted by $C_n$ and $S_{n-1}$ respectively.
We have the initial conditions $\tilde{T}_0(x)=2$, $\tilde{T}_1(x)=x$
and $\tilde{U}_0(x)=0$, $\tilde{U}_1(x)=1$, together with the recurrence relations
\begin{align*}
\tilde{T}_{n+2}(x)=x\tilde{T}_{n+1}(x)-\tilde{T}_n(x),\quad
\tilde{U}_{n+2}(x)=x\tilde{U}_{n+1}(x)-\tilde{U}_n(x).
\end{align*}
The polynomials $\tilde{T}_n$ and $\tilde{U}_n$ obey the following identities:
$$
\begin{array}{r@{\hspace{0.3em}}lr@{\hspace{0.3em}}l}
\tilde{T}_n(2\cos\theta)  &=2\cos(n\theta),&
\tilde{U}_n(2\cos\theta)  &\displaystyle =\frac{\sin(n\theta)}{\sin\theta},
\\[2ex]
\tilde{T}_n(2\cosh\theta) &=2\cosh(n\theta),&
\tilde{U}_n(2\cosh\theta) &\displaystyle =\frac{\sinh(n\theta)}{\sinh\theta},
\\[2ex]
\tilde{T}_n(-2\cosh\theta) &=2(-1)^n\cosh(n\theta),&
\tilde{U}_n(-2\cosh\theta) &\displaystyle =(-1)^{n-1}\frac{\sinh(n\theta)}{\sinh\theta},
\end{array}
$$
and $\tilde{U}_n$ admits the following explicit expression, for any $x\in\mathbb{C}$,
\begin{align*}
\tilde{U}_n(x) & =\sum_{k=0}^{\lfloor (n-1)/2\rfloor} (-1)^k \binom{n-k-1}{k} x^{n-2k-1}.
\end{align*}
As a byproduct, the sequence $(u_n)_{n\in\mathbb{N}}$ introduced in
Theorem~\ref{theorem-gene} is directly related to the family of polynomials
$(\tilde{U}_n)_{n\in\mathbb{N}}$: in fact, $u_n=\tilde{U}_n(a)$. Then
$A=\tilde{U}_{m+1}(a)+\tilde{U}_m(a)$ with $m=(n-1)/2$. This formula
holds for any $a\in\mathbb{C}$. The complement of Theorem~\ref{theorem-gene} is
the following statement.
%
\begin{theorem}
Fix a complex number $a$ and an odd positive integer $n$. The
factor $A$ in the relationship $S_n=Ax_{(n-1)/2}$ related to any sequence
$(x_n)_{n\in\mathbb{N}}$ satisfying the recurrence relation
$x_{n+2}=ax_{n+1}-x_n$ is explicitly given by
\begin{align*}
A =\!\!\sum_{k=0}^{\lfloor (n-1)/4\rfloor} \!\!(-1)^k \left[\binom{(n-3)/2-k}{k} a^{(n-3)/2-2k}
+\binom{(n-1)/2-k}{k} a^{(n-1)/2-2k}\right]\!.
\end{align*}
\end{theorem}
%
We then retrieve~(\ref{proportionA3}). We point out that when $a$ is an
integer, so is $A$.
The reader can compare this expression of $A$ with the various trigonometric and
hyperbolic representations obtained in Section~\ref{subsection-solving-cont}.

Actually, the relationship $S_n=Ax_m$ is nothing but the translation
of a well-known identity concerning the sum of trigonometric (or hyperbolic)
functions. Indeed, recall that, e.g.,
\begin{align*}
\sum_{k=0}^n \cos(k\theta)
=\frac{\sin((n+1)\theta/2)}{\sin(\theta/2)}\,\cos(n\theta/2),
\quad
\sum_{k=0}^n \sin(k\theta)
=\frac{\sin((n+1)\theta/2)}{\sin(\theta/2)}\,\sin(n\theta/2).
\end{align*}
By means of Cheybyshev polynomials, for $n=2m+1$, these two formulas
can be rewritten as
\begin{align}
\sum_{k=0}^{n-1} \tilde{T}_k(2\cos\theta)
& =\frac{\sin((m+1/2)\theta)}{\sin(\theta/2)}\,2\cos(m\theta)
=\frac{\sin(m\theta)+\sin((m+1)\theta)}{\sin\theta}\,2\cos(m\theta)
\nonumber\\
& =[\tilde{U}_m(2\cos\theta)+\tilde{U}_{m+1}(2\cos\theta)]\,\tilde{T}_m(2\cos\theta),
\label{cheb1}\\
\sum_{k=0}^{n-1} \tilde{U}_k(2\cos\theta)
& =\frac{\sin((m+1/2)\theta)}{\sin(\theta/2)}\,\frac{\sin(m\theta)}{\sin\theta}
=\frac{\sin(m\theta)+\sin((m+1)\theta)}{\sin\theta}\,\frac{\sin(m\theta)}{\sin\theta}
\nonumber\\
& =[\tilde{U}_m(2\cos\theta)+\tilde{U}_{m+1}(2\cos\theta)]\,\tilde{U}_m(2\cos\theta).
\label{cheb2}
\end{align}
Formulas~(\ref{cheb1}) and~(\ref{cheb2}) can be simply written,
with $A=\tilde{U}_{m+1}(a)+\tilde{U}_m(a)$, as
\begin{align}\label{proportion-cheb}
\sum_{k=0}^{n-1} \tilde{T}_k(a)=A \,\tilde{T}_m(a)\;\text{ and }\;
\sum_{k=0}^{n-1} \tilde{U}_k(a)=A \,\tilde{U}_m(a).
\end{align}
Since any sequence $(x_n)_{n\in\mathbb{N}}$ satisfying the recurrence
$x_{n+2}=ax_{n+1}-x_n$
is a linear combination of the sequences $(\tilde{T}_n(a))_{n\in\mathbb{N}}$ and
$(\tilde{U}_n(a))_{n\in\mathbb{N}}$, the relationships~(\ref{proportion-cheb}) can
be extended to the sequence $(x_n)_{n\in\mathbb{N}}$, that is
\begin{align*}
\sum_{k=0}^{n-1} x_k=Ax_m.
\end{align*}

\subsection{Link with Fibonacci and Lucas polynomials}\label{subsection-fibo-lucas}

We recall the definition of the families of Fibonacci polynomials
$(F_n)_{n\in\mathbb{N}}$ and Lucas  polynomials $(L_n)_{n\in\mathbb{N}}$
through the expansion of the hyperbolic quantities $\cosh(n\theta)$
and $\sinh(n\theta)$ as polynomials of $\sinh\theta$ (see \cite{lucaspol},
\cite{wikifibopol}):
\begin{align*}
F_n(2\sinh\theta) =\begin{cases}
\displaystyle\frac{\sinh(n\theta)}{\cosh\theta} & \text{if $n$ is even,}
\\[2ex]
\displaystyle\frac{\cosh(n\theta)}{\cosh\theta} & \text{if $n$ is odd,}
\end{cases}
\qquad
L_n(2\sinh\theta) =\begin{cases}2\cosh(n\theta) & \text{if $n$ is even,}
\\[.5ex]
2\sinh(n\theta) & \text{if $n$ is odd.}
\end{cases}
\end{align*}
They are characterized by the initial conditions $F_0(x)=0$, $F_1(x)=1$,
\mbox{$L_0(x)=2$,} $L_1(x)=x$ and the relationships, for any $n\in\mathbb{N}$,
\begin{align*}
F_{n+2}(x)=xF_{n+1}(x)+F_n(x) \;\text{ and }\; L_{n+2}(x)=xL_{n+1}(x)+L_n(x).
\end{align*}
The polynomials $F_n$ and $L_n$ admit the following explicit expansions:
\begin{align*}
F_n(x) =\sum_{k=0}^{\lfloor (n-1)/2\rfloor} \binom{n-k-1}{k} x^{n-2k-1},
\qquad
L_n(x) =\sum_{k=0}^{\lfloor n/2\rfloor} \frac{n}{n-k}\binom{n-k}{k} x^{n-2k}.
\end{align*}
The sequences $(u_n)_{n\in\mathbb{N}}$ and $(v_n)_{n\in\mathbb{N}}$ introduced
in Theorem~\ref{theorem-genebis}
are directly related to the families of polynomials
$(F_n)_{n\in\mathbb{N}}$ and $(L_n)_{n\in\mathbb{N}}$: indeed, $u_n=F_n(a)$ and
$v_n=L_n(a)$. The two expressions of the coefficient $A$ we have obtained
in~(\ref{proportionA2ter}) and~(\ref{proportionA2quater}) can be rewritten as
\begin{align*}
A=F_{(n+p-1)/2}(a)+F_{(n-p+1)/2}(a)
=\begin{cases}
L_{(p-1)/2}(a)\,F_{n/2}(a) & \text{if $n$ is a multiple of $4$,}
\\[.2ex]
F_{(p-1)/2}(a)\,L_{n/2}(a) & \text{if $n$ is a multiple of $4$ plus $2$.}
\end{cases}
\end{align*}
We recall that here $n$ is an even integer and $p$ an odd integer such that
$(n+p-1)$ is a multiple of $4$. Hence $n/2$ and $(p-1)/2$ have the
same parity. The first above equality explicitly yields
\begin{align*}
A=\sum_{k=0}^{\lfloor (n+p-3)/4\rfloor} \binom{(n+p-3)/2-k}{k} a^{(n+p-3)/2-2k}
+\sum_{k=0}^{\lfloor (n-p-1)/4\rfloor} \binom{(n-p-1)/2-k}{k} a^{(n-p-1)/2-2k}.
\end{align*}
This formula holds only for those $a\in\mathbb{C}$ which can be written
as $a=\rho-1/\rho$ where $\rho$ satisfies the equation
$\rho^{p+1}-\rho^p-\rho-1=0$.
%
\begin{remark}
The equality between both expressions~(\ref{proportionA2ter})
and~(\ref{proportionA2quater}) is actually the rephrasing of an elementary identity
concerning the hyperbolic functions. Indeed,
let $s$ and $t$ be two integers with the same parity. Then,
\begin{align*}
F_{s+t}(2\sinh\theta)+F_{s-t}(2\sinh\theta)
& = \frac{\sinh((s+t)\theta)+\sinh((s-t)\theta)}{\cosh\theta}
\\
& = 2\cosh(t\,\theta)\,\frac{\sinh(s\,\theta)}{\cosh\theta}
= 2\sinh(s\,\theta)\,\frac{\cosh(t\,\theta)}{\cosh\theta}
\\
& = \begin{cases}
L_t(2\sinh\theta)F_s(2\sinh\theta) & \text{if $s$ and $t$ are even,}
\\
L_s(2\sinh\theta)F_t(2\sinh\theta) & \text{if $s$ and $t$ are odd.}
\end{cases}
\end{align*}
More concisely,
\begin{align*}
F_{s+t}(a)+F_{s-t}(a)=\begin{cases}
L_t(a)\,F_s(a) & \text{if $s$ and $t$ are even,}
\\
F_t(a)\,L_s(a) & \text{if $s$ and $t$ are odd.}
\end{cases}
\end{align*}
Choosing $s=n/2$ and $t=(p-1)/2$ in the foregoing equality confirms
the equality of~(\ref{proportionA2ter}) and~(\ref{proportionA2quater}).
\end{remark}
%
\begin{remark}
It is well-known that the Fibonacci and Lucas polynomials are
simply related to the Chebyshev polynomials by
\begin{align*}
\tilde{T}_n(ix)=i^nL_n(x)\;\text{ and }\:
\tilde{U}_n(ix)=i^{n-1}F_n(x).
\end{align*}
Despite this connection, because of the factor $i^n$,
there is no simple link between, e.g., the sums $\sum_{k=0}^{n-1}\tilde{U}_k$
and $\sum_{k=0}^{n-1}F_k$. Hence, specific studies were necessary for the
cases $b=1$ and $b=-1$.
\end{remark}
%

\section{Conclusion}\label{section-conclusion}

In this article, we have considered sequences $(x_n)_{n\in\mathbb{N}}$
defined by two initial conditions
$x_0$ and $x_1$ and a second order linear recurrence relation:
$x_{n+2}=ax_{n+1}+bx_n$ for any $n\in\mathbb{N}$. In several cases, for any $n\in\mathbb{N}$,
we have been able to get an integer $m$ and a number $A$ such that
$\sum_{k=0}^{n-1}x_k=Ax_m$. More precisely, such relationship is tractable
for any $a\in\mathbb{C}$ when $b=-1$, and for certain specific values of $a$
when $b=1$. The coefficient $A$ can be expressed in terms of Chebyshev
polynomials when $b=1$ and in terms of Fibonacci polynomials when $b=-1$.

In the table below, we sum up the values of $m$ and $A$
obtained throughout this study.

\begin{center}
\begin{tabular}
{|c@{\hspace{.3em}}|@{\hspace{.3em}}c@{\hspace{.3em}}|@{\hspace{.3em}}c@{\hspace{.2em}}|}
\hline
recurrence & integer $m$ & coefficient $A$
\\
\hline\hline
$x_{n+2}=x_{n+1}+x_n$ &
\parbox{.3\linewidth}{\begin{center}\mbox{}\\[-2.5ex]$n/2+1$\\[-.2ex]
\small ($n$ multiple of $4$ plus $2$)\end{center}} & $\mathcal{L}_{n/2}$
\\[-1.5ex]
\hline
\boite{$x_{n+2}=-x_{n+1}+x_n$}  &
\parbox{.3\linewidth}{\begin{center}\mbox{}\\[-2.5ex]$n/2-2$\\[-.2ex]
\small ($n$ multiple of $4$ plus $2$)\end{center}}
& \boite{$\mathcal{L}_{n/2}$}
\\[-1.5ex]
\hline
\boite{$x_{n+2}=2x_{n+1}+x_n$}  & \boite{$n/2$     \;\small ($n$ multiple of $4$)} & \boite{$2\,F_{n/2}(2)$}
\\[-1.5ex]
\hline
\boite{$x_{n+2}=2x_{n+1}-x_n$}  & \boite{$(n-1)/2$ \;\small (odd $n$)}  & \boite{$n$}
\\[-1.5ex]
\hline
\boite{$x_{n+2}=-2x_{n+1}-x_n$} & \boite{$(n-1)/2$ \;\small (odd $n$)}  & \boite{$(-1)^{(n-1)/2}$}
\\[-1.5ex]
\hline
\boite{$x_{n+2}=ax_{n+1}-x_n$}  & \boite{$(n-1)/2$ \;\small (odd $n$)}
& \boitelarge{$\tilde{U}_{(n+1)/2}(a)+\tilde{U}_{(n-1)/2}(a)$}
\\[-1.5ex]
\hline
\parbox{.3\linewidth}{\mbox{}\quad $x_{n+2}=ax_{n+1}+x_n$\\
\small for certain values of $a\in\mathbb{C}$}
&
\boite{$(n+p-1)/2$\\ \small ($n+p-1$ multiple of $4$)}
& \boitelarge{$F_{(n+p-1)/2}(a)$\\[.5ex]$+\,F_{(n-p+1)/2}(a)$}
\\
\hline
\end{tabular}
\end{center}

\appendix
\section{Appendix}

\setcounter{equation}{0}\renewcommand{\theequation}{\thesection.\arabic{equation}}
\setcounter{theorem}{0}\renewcommand{\thetheorem}{\thesection.\arabic{theorem}}
\setcounter{theorem}{0}\renewcommand{\theproposition}{\thesection.\arabic{proposition}}
\setcounter{theorem}{0}\renewcommand{\thelemma}{\thesection.\arabic{lemma}}

\subsection{An example}

In this appendix, we go back to the example of the recurrence
$x_{n+2}=ax_{n+1}+x_n$ with $a=\rho-1/\rho$ where $\rho$ is the positive root of
$r^6-r^5-r-1=0$; see the case $p=5$ of Examples~\ref{examples}.
Starting from $\rho^6=\rho^5+\rho+1$, we get
\begin{align*}
a&=\frac{\rho^2-1}{\rho}=\frac{\rho^3-\rho}{\rho^2},
\\
a^2&=\frac{\rho^4-2\rho^2+1}{\rho^2},
\\
a^3&=\frac{\rho^6-3\rho^4+3\rho^2-1}{\rho^3}=\frac{\rho^4-3\rho^3+3\rho+1}{\rho^2}.
\end{align*}
We see that $a$ solves the following equation of degree $3$:
\begin{align*}
a^3-a^2+3a-2=0.
\end{align*}
Appealing to the famous Cardano's formula, we can extract the exact real value
of $a$:
\begin{align*}
a=\frac{1}{6}\left(\!\sqrt[3]{\raisebox{0ex}[2ex]{$116+12\sqrt{321}$}}
-\sqrt[3]{\raisebox{0ex}[2ex]{$-116+12\sqrt{321}$}}\,\right)
+\frac{1}{3}.
\end{align*}
As an example of the formulas $S_n=(u_{n/2+2}+u_{n/2-2})x_{n/2+2}
=(a^2+2)u_{n/2}x_{n/2+2}$, we examine the case $n=8$. We begin by writing
the first terms of the sequence $(x_n)_{n\in\mathbb{N}}$:
\begin{align*}
x_2&=x_0+ax_1,
\\
x_3&=ax_0+(a^2+1)x_1,
\\
x_4&=(a^2+1)x_0+(a^3+2a)x_1,
\\
x_5&=(a^3+2a)x_0+(a^4+3a^2+1)x_1,
\\
x_6&=(a^4+3a^2+1)x_0+(a^5+4a^3+3a)x_1,
\\
x_7&=(a^5+4a^3+3a)x_0+(a^6+5a^4+6a^2+1)x_1.
\end{align*}
For $x_0=0$ and $x_1=1$, we derive $u_2=a$, $u_4=a^3+2a$ and
$u_6=a^5+4a^3+3a$. We observe that $u_2+u_6=a^5+4a^3+4a=(a^2+2)u_4$.
The sum of all numbers from $x_0$ to $x_7$ is given by
\begin{align}\label{equality1}
S_8=(a^5+a^4+5a^3+4a^2+6a+4)x_0+(a^6+a^5+6a^4+5a^3+10a^2+6a+4)x_1.
\end{align}
On the other hand,
\begin{align}\label{equality2}
(u_2+u_6)x_6&=(a^5+4a^3+4a)[(a^4+3a^2+1)x_0+(a^5+4a^3+3a)x_1]
\nonumber\\
&=(a^9+7a^7+17a^5+16a^3+4a)x_0+(a^{10}+8a^8+23a^6+28a^4+12a^2)x_1.
\end{align}
Our aim is to check that (\ref{equality1}) and (\ref{equality2}) coincide.
Using $a^3=a^2-3a+2$, we successively obtain
\begin{align*}
a^4    &=a^3-3a^2+2a  =-2a^2-a+2,
\\
a^5    &=-2a^3-a^2+2a =-3a^2+8a-4,
\\
a^6    &=-3a^3+8a^2-4a=5a^2+5a-6,
\\
a^7    &=5a^3+5a^2-6a =10a^2-21a+10,
\\
a^8    &=10a^3-21a^2+10a=-11a^2-20a+20,
\\
a^9    &=-11a^3-20a^2+20a=-31a^2+53a-22,
\\
a^{10} &=-31a^3+53a^2-22a=22a^2+71a-62.
\end{align*}
We easily get
\begin{align*}
u_2+u_6&=a^5+4a^3+4a=a^2+4,
\\
a^5+a^4+5a^3+4a^2+6a+4&=a^9+7a^7+17a^5+16a^3+4a=4a^2-2a+12,
\\
a^6+a^5+6a^4+5a^3+10a^2+6a+4&=a^{10}+8a^8+23a^6+28a^4+12a^2=5a^2-2a+16.
\end{align*}
As a result, we just checked that
\begin{align*}
S_8=(u_2+u_6)x_6=(4a^2-2a+12)x_0+(5a^2-2a+16)x_1
\end{align*}
and we simply have the relationship
\begin{align*}
S_8=(a^2+4)x_6.
\end{align*}

\subsection{An algebraic equation}

In this second appendix, we still consider the example of the recurrence relation
$x_{n+2}=ax_{n+1}+x_n$ treated in Theorem~\ref{theorem-genebis}. Therein,
$a=\rho-1/\rho$ where $\rho$ is a root of $r^{p+1}-r^p-r-1=0$,
$p$ being an odd integer. Of course, since $\rho$ is an algebraic number,
so is $a$. Our aim is to write out an algebraic equation satisfied by $a$.
We distinguish the two cases when $p$ is a multiple of $4$ plus or
minus $1$: $p=4q+1$ or $p=4q-1$ say.

\subsubsection*{\large \sl The case $p=4q+1$}

Suppose that $p$ is of the form $p=4q+1$ for a positive integer $q$.
We have $\rho^{4q+2}-\rho^{4q+1}-\rho-1=0$. The equation
we have obtained for $a$ is an equation of degree $2q+1$.
Our method is heuristic: we search an equation of the form
\begin{align*}
a^{2q+1}=\alpha_{2q}a^{2q}+\alpha_{2q-1}a^{2q-1}+\dots+\alpha_1 a+\alpha_0
=\sum_{j=0}^{2q} \alpha_j a^j
\end{align*}
where $\alpha_0,\alpha_1,\dots,\alpha_{2q}$ are some integers that we are
going to determine.

Let us evaluate the successive powers of $a$: for any positive integer $j$,
\begin{align*}
a^j=\sum_{k=0}^j (-1)^{j+k}\binom{j}{k} \rho^{2k-j}.
\end{align*}
In particular, for $j=2q+1$,
\begin{align*}
a^{2q+1}&=-\sum_{k=0}^{2q+1} (-1)^{k}\binom{2q+1}{k} \rho^{2k-2q-1}
=\rho^{2q+1}-\rho^{-2q-1}-\sum_{k=1}^{2q} (-1)^{k}\binom{2q+1}{k} \rho^{2k-2q-1}.
\end{align*}
Since $\rho^{4q+2}-\rho^{4q+1}-\rho-1=0$, we have
$\rho^{2q+1}-\rho^{-2q-1}=\rho^{2q}+\rho^{-2q}$. Then,
\begin{align*}
a^{2q+1}=\rho^{2q}+\rho^{-2q}-\sum_{k=1}^{2q} (-1)^{k}\binom{2q+1}{k} \rho^{2k-2q-1}.
\end{align*}
The equation we are looking for writes as
\begin{align}\label{sum-inter}
\sum_{j=0}^{2q}\alpha_j\sum_{k=0}^{j} (-1)^{j+k}\binom{j}{k} \rho^{2k-j}
+\sum_{k=1}^{2q} (-1)^{k}\binom{2q+1}{k} \rho^{2k-2q-1}-(\rho^{2q}+\rho^{-2q})=0.
\end{align}
The double sum of the left-hand side of~(\ref{sum-inter}) can be written as follows:
\begin{align}
\sum_{j=0}^{2q}\alpha_j\sum_{k=0}^{j} (-1)^{j+k}\binom{j}{k} \rho^{2k-j}
&=\sum_{i=-2q}^{2q} \Bigg[\sum_{j,k\in\mathbb{N}:\,0\le k\le j\le 2q\atop 2k-j=i}
(-1)^{j+k}\binom{j}{k}\alpha_j\Bigg] \rho^i
\nonumber\\
&=\sum_{i=-2q}^{2q} \Bigg[\sum_{j\in\mathbb{N}:\,|i|\le j\le 2q,\atop  i,j
\text{\,\tiny have the same parity}}
(-1)^{(i+3j)/2}\binom{j}{(i+j)/2}\alpha_j\Bigg] \rho^i
\nonumber\\
&=\alpha_{2q}(\rho^{2q}+\rho^{-2q})
+\sum_{i\in\mathcal{E}:\atop |i|\le 2q-1}
\Bigg[\sum_{j\in\mathcal{E}:\atop |i|\le j\le 2q}
(-1)^{(j-i)/2}\binom{j}{(i+j)/2}\alpha_j\Bigg]\rho^i
\nonumber\\
&\phantom{=\;}+\sum_{i\in\mathcal{O}:\atop |i|\le 2q-1}
\Bigg[\sum_{j\in\mathcal{O}:\atop |i|\le j\le 2q}
(-1)^{(j-i)/2}\binom{j}{(i+j)/2}\alpha_j\Bigg] \rho^i.
\label{sum1-inter}
\end{align}
In the last equality, $\mathcal{E}$ and $\mathcal{O}$ denote
the sets of even and odd integers respectively.
By using symmetry properties of the binomial coefficients, it is
enough in~(\ref{sum1-inter}) to sum over positive integers:
\begin{align}
\sum_{j=0}^{2q}\alpha_j\sum_{k=0}^{j} (-1)^{j+k}\binom{j}{k} \rho^{2k-j}
&=\alpha_{2q}(\rho^{2q}+\rho^{-2q})
+\sum_{j\in\mathcal{E}:\atop 0\le j\le 2q}
(-1)^{j/2}\binom{j}{j/2}\alpha_j
\nonumber\\
&\phantom{=\;}+\sum_{i\in\mathcal{E}:\atop 1\le i\le 2q-1}
\Bigg[\sum_{j\in\mathcal{E}:\atop i\le j\le 2q}
(-1)^{(j-i)/2}\binom{j}{(i+j)/2}\alpha_j\Bigg]
(\rho^i+\rho^{-i})
\nonumber\\
&\phantom{=\;}+\sum_{i\in\mathcal{O}:\atop 1\le i\le 2q-1}
\Bigg[\sum_{j\in\mathcal{O}:\atop i\le j\le 2q}
(-1)^{(j-i)/2}\binom{j}{(i+j)/2}\alpha_j\Bigg]
(\rho^i-\rho^{-i}).
\label{sum1bis-inter}
\end{align}
On the other hand,
\begin{align}\label{sum2-inter}
\sum_{k=1}^{2q} (-1)^{k}\binom{2q+1}{k} \rho^{2k-2q-1}
=\sum_{i\in\mathcal{O}:\atop 1\le i\le 2q-1}
(-1)^{(i+1)/2+q}\binom{2q+1}{(i+1)/2+q} (\rho^i-\rho^{-i}).
\end{align}
Plugging~(\ref{sum1bis-inter}) and~(\ref{sum2-inter}) into~(\ref{sum-inter})
yields
\begin{align}
\noalign{\noindent$\displaystyle
(\alpha_{2q}-1)(\rho^{2q}+\rho^{-2q})
+\sum_{j\in\mathcal{E}:\atop 0\le j\le 2q}
(-1)^{j/2}\binom{j}{j/2}\alpha_j
+\sum_{i\in\mathcal{E}:\atop 1\le i\le 2q-1}
\Bigg[\sum_{j\in\mathcal{E}:\atop i\le j\le 2q}
(-1)^{(j-i)/2}\binom{j}{(i+j)/2}\alpha_j\Bigg]
(\rho^i+\rho^{-i})
$}
\noalign{\noindent$\displaystyle
+\sum_{i\in\mathcal{O}:\atop 1\le i\le 2q-1}
\Bigg[(-1)^{(i+1)/2+q}\binom{2q+1}{(i+1)/2+q}
+\sum_{j\in\mathcal{O}:\atop i\le j\le 2q}
(-1)^{(j-i)/2}\binom{j}{(i+j)/2}\alpha_j\Bigg]
(\rho^i-\rho^{-i})=0.
$}
\label{sum3-inter}
\end{align}
With obvious notations, Equation~(\ref{sum3-inter}) is clearly of the form
\begin{align*}
\sum_{i=0}^{2q}\beta_i(\rho^i+(-1)^i\rho^{-i})=0
\end{align*}
which is satisfies as soon as the $\beta_i$'s vanish (this is a
sufficient condition which is not \textit{a priori} necessary).
Hence, (\ref{sum3-inter}) is satisfied when $\alpha_{2q}=1$ and
\begin{align*}
\begin{cases}
\displaystyle
\sum_{j\in\mathcal{E}:\atop i\le j\le 2q}
(-1)^{(j-i)/2}\binom{j}{(i+j)/2}\alpha_j
&
\!\!\!\displaystyle
=0\quad\text{for $i\in\{0,2,4,\dots,2q-2\}$,}
\\[3ex]
\displaystyle
\sum_{j\in\mathcal{O}:\atop i\le j\le 2q}
(-1)^{(j-i)/2}\binom{j}{(i+j)/2}\alpha_j
&
\!\!\!\displaystyle
=(-1)^{(i-1)/2+q}\binom{2q+1}{(i+1)/2+q}
\quad\text{for $i\in\{1,3,5,\dots,2q-1\}$.}
\end{cases}
\end{align*}
We rewrite this system by changing the indices $i$ and $j$
into $2i$ and $2j$ in the first equation, $(2i+1)$ and $(2j+1)$
in the second one:
\begin{align}
\begin{cases}
\displaystyle
\sum_{j=i}^{q-1} (-1)^j\binom{2j}{j-i}\alpha_{2j}
&
\!\!\!\displaystyle
=(-1)^{q+1}\binom{2q}{q-i}
\quad\text{for $i\in\{0,1,2,\dots,q-1\}$,}
\\[3ex]
\displaystyle
\sum_{j=i}^{q-1} (-1)^j \binom{2j+1}{j-i}\alpha_{2j+1}
&
\!\!\!\displaystyle
=(-1)^q\binom{2q+1}{q-i}
\,\,\quad\text{for $i\in\{0,1,2,\dots,q-1\}$.}
\end{cases}
\nonumber\\[-10ex]
\label{system}
\\[3ex]
\nonumber
\end{align}
Let us introduce the principal matrices related to this system:
\begin{align*}
T_1=\left[(-1)^j\binom{2j}{j-i}\right]_{0\le i,j\le q-1}
\;\text{ and }\;
T_2=\left[(-1)^j\binom{2j+1}{j-i}\right]_{0\le i,j\le q-1}
\end{align*}
as well as the matrices of the unknowns and those of the right-hand sides
of~(\ref{system}):
\begin{align*}
A_1=\left[\alpha_{2i}\right]_{0\le i\le q-1},\quad
A_2=\left[\alpha_{2i+1}\right]_{0\le i\le q-1},
\end{align*}
\begin{align*}
B_1=(-1)^{q+1}\left[\binom{2q}{q-i}\right]_{0\le i\le q-1},\quad
B_2=(-1)^q\left[\binom{2q+1}{q-i}\right]_{0\le i\le q-1}.
\end{align*}
Recalling the convention that $\binom{\ell}{k}=0$ if $k>\ell$,
we notice that the matrices $T_1$ and $T_2$ are triangular.
With these settings at hand, the system~(\ref{system}) simply reads
\begin{align*}
T_1A_1=B_1 \quad\text{and}\quad T_2A_2=B_2.
\end{align*}
We need to compute the inverse of $T_1$ and $T_2$.
By performing empirical computations with the help of Maple, we found
that the entries of $T_1^{-1}$ and $A_1=T_1^{-1}B_1$ are
respectively
\begin{align*}
(-1)^j\binom{2j-2i}{j-i}\binom{2j}{2i}/\binom{2j-1}{j-i}
\quad\text{and}\quad\binom{2q}{2i}\binom{2q-2i}{q-i}/\binom{2q-1}{q-i}.
\end{align*}
We found analogous expressions $T_2^{-1}$ and $A_2=T_2^{-1}B_2$.
All of them can be simplified. We state the results in the
Propositions~\ref{prop1} and~\ref{prop2}.
%
\begin{proposition}\label{prop1}
We have
\begin{align*}
T_1^{-1} &=\left[(-1)^j\frac{2j}{i+j}
\binom{i+j}{2i}\right]_{0\le i,j\le q-1}\!,
\\[1ex]
T_2^{-1} &=\left[(-1)^j\frac{2j+1}{i+j+1}
\binom{i+j+1}{2i+1}\right]_{0\le i,j\le q-1}\!.
\end{align*}
In $T_1^{-1}$ above, we adopt the convention that $(2j)/(i+j)=1$ if $i=j=0$.
\end{proposition}
%
\proof

We only prove the formula related to $T_1^{-1}$, the proof of that
related to $T_2^{-1}$ being quite similar.
We check that the product of upper triangular matrices
\begin{align}\label{product}
\left[(-1)^j\binom{2j}{i+j}\right]_{0\le i,j\le q-1}\times
\left[(-1)^j\frac{2j}{i+j} \binom{i+j}{2i}\right]_{0\le i,j\le q-1}
\end{align}
exactly coincides with the unit matrix of type $q\times q$,
that is, $T_1T_1^{-1}=I_q$. For this, we have to compute
its upper entries: for $0\le i\le q-1$ and $1\le j\le q-1$ such that
$i\le j$,
\begin{align*}
\sum_{k=0}^{q-1} (-1)^{j+k} \frac{2j}{k+j} \binom{2k}{k+i} \binom{k+j}{2k}
&=2j\sum_{k=i}^j (-1)^{j+k} \frac{(k+j-1)!}{(j-k)!(k-i)!(k+i)!}
\\
&=(-1)^j \frac{2j}{(j-i)!}\sum_{k=i}^j (-1)^k \frac{(k+j-1)!}{(k+i)!}\binom{j-i}{k-i}
\\
&=(-1)^{i+j} \frac{2j}{(j-i)!}\sum_{k=0}^{j-i}
(-1)^k \frac{(k+i+j-1)!}{(k+2i)!}\binom{j-i}{k}.
\end{align*}
In order to evaluate the last sum, we use the lemma below:
applying~(\ref{sum-binom}) to $\ell=2i+j$ and $m=j-i$ immediately yields,
for $0\le i\le q-1$ and $1\le j\le q-1$ such that $i\le j$,
\begin{align*}
\sum_{k=0}^{q-1} (-1)^{j+k} \frac{2j}{k+j} \binom{2k}{k+i} \binom{k+j}{2k}=0.
\end{align*}
With the convention $(2j)/(k+j)=1$ for $j=k=0$, the foregoing sum is equal
to $1$. This proves that the product~(\ref{product}) coincides with the
unit matrix.
\qed
%
\begin{lemma}
For any positive integers $\ell$ and $m$ such that $\ell\ge m$,
we have
\begin{align}\label{sum-binom}
\sum_{k=0}^m (-1)^k \frac{(k+\ell-1)!}{(k+\ell-m)!}\binom{m}{k}=0.
\end{align}
\end{lemma}
%
\proof

By observing that $\displaystyle\frac{\mathrm{d}^{m-1}}{\mathrm{d}^{m-1}x} (x^{k+\ell-1})
=\frac{(k+\ell-1)!}{(k+\ell-m)!} \,x^{k+\ell-m}$, we get
\begin{align*}
\sum_{k=0}^m (-1)^k \frac{(k+\ell-1)!}{(k+\ell-m)!}\binom{m}{k} x^k
=x^{m-\ell}\frac{\mathrm{d}^{m-1}}{\mathrm{d}^{m-1}x}
\left(\sum_{k=0}^m (-1)^k\binom{m}{k} x^{k+\ell-1}\right)
=x^{m-\ell} \frac{\mathrm{d}^{m-1}}{\mathrm{d}^{m-1}x}[(1-x)^mx^{\ell-1}].
\end{align*}
It is clear that $1$ is a root of $(1-x)^mx^{\ell-1}$ of order $m$, so that
the $(m-1)$th derivative of $(1-x)^mx^{\ell-1}$ vanishes at $x=1$.
This finishes the proof of~(\ref{sum-binom}).
\qed
%
\begin{proposition}\label{prop2}
We have
\begin{align*}
A_1=\left[\frac{2q}{q+i} \binom{q+i}{2i}\right]_{0\le i\le q-1}\!,
\quad A_2 =-\left[\frac{2q+1}{q+i+1} \binom{q+i+1}{2i+1}\right]_{0\le i\le q-1}\!.
\end{align*}
\end{proposition}
%
\proof

We only prove the formula related to $A_1$, the proof of that
related to $A_2$ being quite similar.
We have to compute the product $A_1=T_1^{-1}B_1$:
\begin{align*}
\left[(-1)^j\frac{2j}{i+j} \binom{i+j}{2i}\right]_{0\le i,j\le q-1}
\times\left[(-1)^{q+1}\binom{2q}{q-i}\right]_{0\le i\le q-1}.
\end{align*}
The entries of this matrix are
\begin{align*}
\alpha_{2i}=\sum_{k=i}^{q-1} (-1)^{k+q+1} \frac{2k}{k+i} \binom{k+i}{2i} \binom{2q}{q-k}.
\end{align*}
We previously proved that the product~(\ref{product}) coincides with the unit matrix,
that is $T_1T_1^{-1}=I_q$. Of course, in this case, the product is
commutative and we have also $T_1^{-1}T_1=I_q$. This remark entails the
following identity: for any integers $i,j$ such that $0\le i\le j\le q-1$,
\begin{align*}
\sum_{k=i}^j (-1)^{j+k} \frac{2k}{k+i} \binom{k+i}{2i} \binom{2j}{k+j}=\delta_{ij}.
\end{align*}
Noticing that this formula actually does not depend on the number $q$,
it holds true also for $j=q$, which gives, for $0\le i\le q-1$,
\begin{align}\label{sum4}
\sum_{k=i}^q (-1)^k \frac{2k}{k+i} \binom{k+i}{2i} \binom{2q}{k+q}=0.
\end{align}
As a result, by isolating the last term of the previous sum, we finally get,
for $0\le i\le q-1$,
\begin{align*}
\sum_{k=i}^{q-1} (-1)^{k+q+1} \frac{2k}{k+i} \binom{k+i}{2i} \binom{2q}{q-k}
=\frac{2q}{q+i} \binom{q+i}{2i}.
\end{align*}
This proves the formula obtained for $A_1$.
\qed

The conclusion of this analysis is the following equation satisfied by $a$.
%
\begin{theorem}
Let $q$ be a positive integer and $\rho$ be a root of the equation
$r^{4q+2}=r^{4q+1}+r+1$.
The parameter $a=\rho-1/\rho$ is an algebraic number; it solves the equation
\begin{align*}
a^{2q+1} =\sum_{i=0}^{q} \frac{2q}{q+i} \binom{q+i}{2i}a^{2i}
-\sum_{i=0}^{q-1} \frac{2q+1}{q+i+1} \binom{q+i+1}{2i+1}a^{2i+1}.
\end{align*}
\end{theorem}
%
\begin{examples}
For $q\in\{1,2,3,4\}$ (\/$p\in\{5,9,13,17\}$), the corresponding equation
satisfied by $a$ reads
\begin{align*}
&a^3-a^2+3a-2=0,
\\
&a^5-a^4+5a^3-4a^2+5a-2=0,
\\
&a^7-a^6+7a^5-6a^4+14a^3-9a^2+7a-2=0,
\\
&a^9-a^8+9a^7-8a^6+27a^5-20a^4+30a^3-16a^2+9a-2=0.
\end{align*}
\end{examples}

\subsubsection*{\large \sl The case $p=4q-1$}

Assume now that $p=4q-1$ for a positive integer $q$.
We have in this case $\rho^{4q}-\rho^{4q-1}-\rho-1=0$.
The previous calculations cannot be carried out exactly in the same way
and we must work with a slightly modified equation.
In this situation, we observe that $i$ and $-i$ are particular solutions;
we also choose $a$ different from $\pm 2i$ as required in
Theorem~\ref{theorem-genebis}. Thus, the polynomial $r^{4q}-r^{4q-1}-r-1$
can be divided by $r^2+1$ and the quotient is
\begin{align*}
r^{4q-2}-r^{4q-3}-r^{4q-4}+r^{4q-5}+r^{4q-6}-\dots+r^3+r^2-r-1
=r^{4q-2}-\sum_{k=0}^{2q-2} (-1)^k(r^{2k}+r^{2k+1}).
\end{align*}
The equation we have obtained for $a$ is an equation of degree $2q-1$.
As previously, we look for an equation of the form
\begin{align*}
a^{2q-1}=\alpha_{2q-2}a^{2q-2}+\alpha_{2q-3}a^{2q-3}+\dots+\alpha_1 a+\alpha_0
=\sum_{j=0}^{2q-2} \alpha_j a^j
\end{align*}
where $\alpha_0,\alpha_1,\dots,\alpha_{2q-2}$ are some integers that we are
going to determine. The successive powers of $a$ write, for any positive integer $j$,
as
\begin{align*}
a^j=\sum_{k=0}^j (-1)^{j+k}\binom{j}{k} \rho^{2k-j}.
\end{align*}
In particular, for $j=2q-1$,
\begin{align*}
a^{2q-1}&=\sum_{k=0}^{2q-1} (-1)^{k-1}\binom{2q-1}{k} \rho^{2k-2q+1}
=\rho^{2q-1}-\sum_{k=0}^{2q-2} (-1)^{k}\binom{2q-1}{k} \rho^{2k-2q+1}.
\end{align*}
Since $\rho^{4q-2}-\sum_{k=0}^{2q-2} (-1)^k(\rho^{2k}+\rho^{2k+1})=0$, we have
\begin{align*}
\rho^{2q-1}=\sum_{k=0}^{2q-2} (-1)^k \rho^{2k-2q+1}+
\sum_{k=1}^{2q-1} (-1)^{k-1} \rho^{2k-2q}.
\end{align*}
Therefore,
\begin{align*}
a^{2q-1}=\sum_{k=1}^{2q-2} (-1)^{k-1}\left[\binom{2q-1}{k}-1\right]\rho^{2k-2q+1}
+\sum_{k=1}^{2q-1} (-1)^{k-1} \rho^{2k-2q}.
\end{align*}
The equation we are looking for becomes
\begin{align}\label{sumbis-inter}
\sum_{j=0}^{2q-2}\alpha_j\sum_{k=0}^{j} (-1)^{j+k}\binom{j}{k} \rho^{2k-j}
+\sum_{k=1}^{2q-2} (-1)^k\left[\binom{2q-1}{k}-1\right]\rho^{2k-2q+1}
+\sum_{k=1}^{2q-1} (-1)^k\rho^{2k-2q}=0.
\end{align}
This is the analogue of~(\ref{sum1bis-inter}).
The double sum in the left-hand side is given by~(\ref{sum1-inter})
upon changing $q$ into $q-1$:
\begin{align}
\sum_{j=0}^{2q-2}\alpha_j\sum_{k=0}^{j} (-1)^{j+k}\binom{j}{k} \rho^{2k-j}
&=\sum_{j\in\mathcal{E}:\atop 0\le j\le 2q-2} (-1)^{j/2}\binom{j}{j/2}\alpha_j
\nonumber\\
&\phantom{==}+\sum_{i\in\mathcal{E}:\atop 1\le i\le 2q-2}
\Bigg[\sum_{j\in\mathcal{E}:\atop i\le j\le 2q-2}
(-1)^{(j-i)/2}\binom{j}{(i+j)/2}\alpha_j\Bigg](\rho^i+\rho^{-i})
\nonumber\\
&\phantom{==}+\sum_{i\in\mathcal{O}:\atop 1\le i\le 2q-3}
\Bigg[\sum_{j\in\mathcal{O}:\atop i\le j\le 2q-3}
(-1)^{(j-i)/2}\binom{j}{(i+j)/2}\alpha_j\Bigg] (\rho^i-\rho^{-i}).
\label{sum1ter-inter}
\end{align}
On the other hand,
\begin{align}
\sum_{k=1}^{2q-1} (-1)^k\rho^{2k-2q}=(-1)^q+
\sum_{i\in\mathcal{E}:\atop 1\le i\le 2q-2}(-1)^{i/2+q} (\rho^i+\rho^{-i}),
\label{sum2bis-inter}
\end{align}
and
\begin{align}
\sum_{k=1}^{2q-2} (-1)^k\left[\binom{2q-1}{k}-1\right]\rho^{2k-2q+1}
=\sum_{i\in\mathcal{O}:\atop 1\le i\le 2q-3}
(-1)^{(i-1)/2+q}\left[\binom{2q-1}{(i-1)/2+q}-1\right] (\rho^i-\rho^{-i}).
\label{sum3bis-inter}
\end{align}
Plugging~(\ref{sum1ter-inter}), (\ref{sum2bis-inter})
and~(\ref{sum3bis-inter}) into~(\ref{sumbis-inter}) yields
\begin{align*}
\noalign{\noindent$\displaystyle
\Bigg[\sum_{j\in\mathcal{E}:\atop 0\le j\le 2q-2}
(-1)^{j/2}\binom{j}{j/2}\alpha_j+(-1)^q\Bigg]
+\sum_{i\in\mathcal{E}:\atop 1\le i\le 2q-2}
\Bigg[\sum_{j\in\mathcal{E}:\atop i\le j\le 2q-2}
(-1)^{(j-i)/2}\binom{j}{(i+j)/2}\alpha_j+(-1)^{i/2+q}\Bigg]
(\rho^i+\rho^{-i})
$}
\noalign{\noindent$\displaystyle
+\sum_{i\in\mathcal{O}:\atop 1\le i\le 2q-3}
\Bigg[\sum_{j\in\mathcal{O}:\atop i\le j\le 2q-3}
(-1)^{(j-i)/2}\binom{j}{(i+j)/2}\alpha_j
+(-1)^{(i-1)/2+q}\left[\binom{2q-1}{(i-1)/2+q}-1\right]
\Bigg](\rho^i-\rho^{-i})=0.
$}
\end{align*}
This equation is of the form
\begin{align*}
\sum_{i=0}^{2q-2}\beta_i(\rho^i+(-1)^i\rho^{-i})=0
\end{align*}
which is satisfied as soon as the $\beta_i$'s vanish, that is
when the $\alpha_j$'s satisfy
\begin{align*}
\begin{cases}
\displaystyle
\sum_{j\in\mathcal{E}:\atop i\le j\le 2q-2}(-1)^{j/2}\binom{j}{(i+j)/2}\alpha_j
&
\!\!\!\!=(-1)^{q-1}\quad\text{for $i\in\{0,2,4,\dots,2q-2\}$,}
\\[3ex]
\displaystyle
\sum_{j\in\mathcal{O}:\atop i\le j\le 2q-3}(-1)^{(j+1)/2}\binom{j}{(i+j)/2}\alpha_j
&
\!\!\!\!\displaystyle=(-1)^{q}\left[\binom{2q-1}{(i-1)/2+q}-1\right]
\quad\text{for $i\in\{1,3,5,\dots,2q-3\}$.}
\end{cases}
\end{align*}
With an obvious change of indices, we rewrite this system as
\begin{align}
\begin{cases}
\displaystyle
\sum_{j=i}^{q-1} (-1)^j\binom{2j}{j-i}\alpha_{2j}
&
\!\!\!\!\displaystyle
=(-1)^{q-1}\quad\text{for $i\in\{0,1,2,\dots,q-1\}$,}
\\[3ex]
\displaystyle
\sum_{j=i}^{q-1} (-1)^j \binom{2j-1}{j-i}\alpha_{2j-1}
&
\!\!\!\!\displaystyle
=(-1)^q\!\left[\binom{2q-1}{q-i}-1\right]\quad\text{for $i\in\{1,2,3,\dots,q-1\}$.}
\end{cases}
\nonumber\\[-10ex]
\label{systembis}
\\[3ex]
\nonumber
\end{align}
By introducing the matrices
\begin{align*}
&T_1=\left[(-1)^j\binom{2j}{j-i}\right]_{0\le i,j\le q-1},\quad
T_2=\left[(-1)^j\binom{2j-1}{j-i}\right]_{1\le i,j\le q-1},
\\
&A_1=\left[\alpha_{2i}\right]_{0\le i\le q-1},\quad
A_2=\left[\alpha_{2i-1}\right]_{1\le i\le q-1},
\\
&B_1=(-1)^{q-1}\left[1\right]_{0\le i\le q-1},\quad
B_2=(-1)^q\left[\binom{2q-1}{q-i}\right]_{1\le i\le q-1}
\!-(-1)^q\left[1\right]_{1\le i\le q-1},
\end{align*}
the system~(\ref{systembis}) can be simply written as
$T_1A_1=B_1$ and $T_2A_2=B_2.$
Observe that the matrix $T_1$ is exactly the same as the $T_1$
in the previous case, while $T_2$ is nearly the same as $T_2$
in the previous case upon changing $q$ into $q-1$.
%
\begin{proposition}\label{prop3}
We have
\begin{align*}
A_1=\left[\binom{q+i-1}{2i}\right]_{0\le i\le q-1}\!,
\quad A_2 =-\left[\binom{q+i-1}{2i-1}\right]_{1\le i\le q-1}\!.
\end{align*}
\end{proposition}
%
\proof
Referring to Proposition~\ref{prop1} for $T_1^{-1}$, we first compute the product $A_1=T_1^{-1}B_1$:
\begin{align*}
\left[(-1)^j\frac{2j}{i+j} \binom{i+j}{2i}\right]_{0\le i,j\le q-1}
\times\left[(-1)^{q-1}\right]_{0\le i\le q-1}.
\end{align*}
The entries of this matrix are
\begin{align}
\alpha_{2i}&=\sum_{k=i}^{q-1} (-1)^{k+q-1} \frac{2k}{k+i} \binom{k+i}{2i}
=\sum_{k=i}^{q-1} (-1)^{k+q-1} \left(2-\frac{2i}{k+i}\right) \binom{k+i}{2i}
\nonumber\\
&=\sum_{k=i}^{q-1} (-1)^{k+q-1} \left[2\binom{k+i}{2i}-\binom{k+i-1}{2i-1}\right]\!.
\label{alpha2i-inter}
\end{align}
By using Pascal's rule, we replace $\binom{k+i}{2i}-\binom{k+i-1}{2i-1}$
by $\binom{k+i-1}{2i}$ which yields
\begin{align*}
\alpha_{2i}&=(-1)^{q-1}\sum_{k=i}^{q-1}
\left[(-1)^{k}\binom{k+i}{2i}-(-1)^{k-1}\binom{k+i-1}{2i}\right]\!.
\end{align*}
Observing that the above sum is
telescopic, we finally obtain the value of $\alpha_{2i}$:
\begin{align}\label{alpha2i}
\alpha_{2i}=\binom{q+i-1}{2i}.
\end{align}
Referring to Proposition~\ref{prop1} for $T_2^{-1}$ (upon changing $q$ into $q-1$),
we now compute the product $A_2=T_2^{-1}B_2$ which can be split into
the difference of two matrices $A_2=A'_2-A''_2$ because of the structure of $B_2$,
where
\begin{align*}
A'_2&=(-1)^{q}\!\left[(-1)^j\frac{2j-1}{i+j-1}\binom{i+j-1}{2i-1}\right]_{1\le i,j\le q-1}
\times\left[\binom{2q-1}{q-i}\right]_{1\le i\le q-1},
\\[2ex]
A''_2&=(-1)^{q}\!\left[(-1)^j\frac{2j-1}{i+j-1}\binom{i+j-1}{2i-1}\right]_{1\le i,j\le q-1}
\times\left[1\right]_{1\le i\le q-1}.
\end{align*}
The entries of $A'_2$ are
\begin{align*}
\alpha'_{2i-1}=\sum_{k=i}^{q-1}(-1)^{k+q}\frac{2k-1}{k+i-1}\binom{k+i-1}{2i-1}\binom{2q-1}{q-k}.
\end{align*}
By expanding the product $T_2^{-1}T_2$ (as in the proof of
Proposition~\ref{prop2}) which coincides with the unit matrix of type
$(q-1)\times (q-1)$ and changing $q$ into $q+1$ therein (since it works for any $q$),
we obtain the analogue of~(\ref{sum4}), namely for any $i$ such that $1\le i\le q-1$,
\begin{align*}
\sum_{k=i}^{q}(-1)^k\frac{2k-1}{k+i-1}\binom{k+i-1}{2i-1}\binom{2q-1}{q-k}=0.
\end{align*}
By isolating the last term in the above sum, we get
\begin{align*}
\sum_{k=i}^{q-1}(-1)^k\frac{2k-1}{k+i-1}\binom{k+i-1}{2i-1}\binom{2q-1}{q-k}
=(-1)^{q-1}\frac{2q-1}{q+i-1}\binom{q+i-1}{2i-1}
\end{align*}
and
\begin{align*}
\alpha'_{2i-1}=-\frac{2q-1}{q+i-1}\binom{q+i-1}{2i-1}.
\end{align*}
On the other hand, the entries of $A''_2$ are
\begin{align*}
\alpha''_{2i-1}=\sum_{k=i}^{q-1}(-1)^{k+q}\frac{2k-1}{k+i-1}\binom{k+i-1}{2i-1}.
\end{align*}
This expression is quite similar to the first sum of~(\ref{alpha2i-inter})
and can be evaluated in the same manner by writing
$\frac{2k-1}{k+i-1}=2-\frac{2i-1}{k+i-1}$ and using Pascal's rule.
We find the analogue of~(\ref{alpha2i}):
\begin{align*}
\alpha''_{2i-1}=-\binom{q+i-2}{2i-1}.
\end{align*}
Finally, the entries of $A_2=A'_2-A''_2$ write as
\begin{align*}
\alpha_{2i-1}=\alpha'_{2i-1}-\alpha''_{2i-1}
=-\frac{2q-1}{q+i-1}\binom{q+i-1}{2i-1}+\binom{q+i-2}{2i-1}
\end{align*}
which can be simplified into
\begin{align*}
\alpha_{2i-1}=-\binom{q+i-1}{2i-1}.
\end{align*}
The proof of Proposition~\ref{prop3} is finished.
\qed

Consequently, we can state the result below.
%
\begin{theorem}
Let $q$ be a positive integer and $\rho$ be a root of
the equation $r^{4q}=r^{4q-1}+r+1$ distinct of $\pm i$.
The parameter $a=\rho-1/\rho$ is an algebraic number; it solves the equation
\begin{align*}
a^{2q-1} =\sum_{i=0}^{q-1} \binom{q+i-1}{2i}a^{2i}
-\sum_{i=1}^{q-1} \binom{q+i-1}{2i-1}a^{2i-1}.
\end{align*}
\end{theorem}
%
\begin{examples}
For $q\in\{1,2,3,4\}$ (\/$p\in\{3,7,11,15\}$), the corresponding equation
satisfied by $a$ reads
\begin{align*}
&a-1=0,
\\
&a^3-a^2+2a-1=0,
\\
&a^5-a^4+4a^3-3a^2+3a-1=0,
\\
&a^7-a^6+6a^5-5a^4+10a^3-6a^2+4a-1=0.
\end{align*}
\end{examples}

\section*{Acknowledgements} I warmly thank my colleague Pierre Schott%
\footnote{\'Ecole Sup\'erieure d'Informatique, d'\'Electronique et d'Automatisme,
P\^ole ARNUM, 9 rue V\'esale, 75005 Paris, France. Email: pierre.schott@esiea.fr}
for having submitted this problem to me and given me many references on the
original trick. He is a physicist in a french engineering school
(College in Informatics, Electronics and  Automatism, Paris).
He is fascinated by magic and he often uses magic for making his lectures (in mathematics,
physics, electronic, computing...) more lively, illustrative and attractive.


\end{document}